\documentclass[review,onefignum,onetabnum]{siamart190516}

\usepackage{lipsum}
\usepackage{amsfonts}
\usepackage{graphicx}
\usepackage{epstopdf}
\usepackage{algorithmic}

\usepackage[rotateright]{rotating}
\usepackage{subfigure,fancybox}
\usepackage{amscd,amstext}
\usepackage{amsmath}
\usepackage{amssymb}
\usepackage{color}
\usepackage{multirow}
\usepackage{mathrsfs}
\usepackage{diagbox}

\usepackage{geometry}
\ifpdf
  \DeclareGraphicsExtensions{.eps,.pdf,.png,.jpg}
\else
  \DeclareGraphicsExtensions{.eps}
\fi

\newsiamremark{rem}{Remark}
\newsiamremark{hypothesis}{Hypothesis}
\crefname{hypothesis}{Hypothesis}{Hypotheses}
\newsiamthm{claim}{Claim}
\newsiamthm{lem}{Lemma}
\newsiamthm{thm}{Theorem}

\headers{splitting algorithm for
optimal control problem}{H. Song, J. Zhang, and Y. Hao}

\title{A splitting algorithm for constrained optimization problems with parabolic equations\thanks{
\funding{The work of H. Song was supported by the NSF of China under the grant
No.11701210, the NSF of Jilin Province under the grants No. 20190103029JH, 20200201269JC,
the education department project of Jilin Province under the grant No. JJKH20211031KJ,
and the fundamental research funds for the Central Universities.
The work of J.C. Zhang was supported by the Natural Science Foundation of Jiangsu Province (Grant BK20210540)
, the Natural Science Foundation of the Jiangsu Higher Education Institutions of China
(No. 21KJB110015, 21KJB110001) and the Startup Foundation for Introducing Talent of NJTech (No. 39804131).
The work of Y.L. Hao was supported
by the NSF of China under the grant No. 11901606.}}}

\author{
Haiming Song\thanks{School of Mathematics, Jilin University, Changchun 130012, China
  (\email{songhaiming@jlu.edu.cn}).}
\and Jiachuan Zhang\thanks{Corresponding author. School of Physical and Mathematical Sciences, Nanjing Tech University, Nanjing 211816, China
  (\email{zhangjc@njtech.edu.cn}).}
\and Yongle Hao\thanks{School of Mathematics and Statistics, Zhoukou Normal University, Zhoukou, 466001, China
  (\email{haoyl@zknu.edu.cn}).}
}

\usepackage{amsopn}

\begin{document}

\maketitle

\begin{abstract}
In this paper, an efficient parallel splitting method is proposed
for the optimal control problem with parabolic equation constraints.
The linear finite element is used to approximate the
state variable and the control variable in spatial direction. And
the Crank-Nicolson scheme is applied to discretize the constraint
equation in temporal direction. For consistency, the trapezoidal
rule and midpoint rule are used to approximate the integrals with
respect to the state variable and the control variable of the
objective function in temporal direction, respectively. Based on
the separable structure of the resulting coupled discretized
optimization system, a full Jacobian decomposition method with
correction is adopted to solve the decoupled subsystems in parallel,
which improves the computational efficiency significantly. Moreover,
the global convergence estimate is established using the discretization
error by the finite element and the iteration error by the full Jacobian
decomposition method with correction. Finally, numerical simulations are
carried out to verify the efficiency of the proposed method.
\end{abstract}



\section{Introduction}\label{INT}

Because of its widespread applications in engineering, mathematical
finance, physics, and life sciences fields, the optimal control
problem with partial differential equation (PDE) constraints has
always been the focus of the scientific computing communities.
Therefore, there exist fruitful research results on this topic
theoretically and numerically \cite{BY17,YWX19},
especially for the elliptic optimal control problem
\cite{GXY17,HPUU10,WAYT21}. Because of the large scale
of discrete system and the limitation of computational resources,
the design of the high accuracy mathematical scheme for parabolic
optimal control problem is difficult, and we refer the readers to
\cite{AF12,HCH11,Kwo17,LW18,LCHH14,LLX17,LMTY04}
and references therein. Based on the finite element approximation,
the Crank-Nicolson scheme, the numerical integration formula, and
the full Jacobian decomposition method with correction \cite{AF12,HHY15},
an efficient numerical algorithm is proposed for the parabolic optimal
control problem in this paper.

Let $y$ and $u$ be the state variable and the control variable,
respectively. Assume $\Omega \subset \mathbb{R}^2$ is a bounded
polygonal domain with Lipschitz boundary and $Q_T=\Omega\times(0,T]$.
Consider the following parabolic optimal control problem
\begin{eqnarray}\label{eq: objection}
\min_{y,u}\mathcal{J}(y,u)=\frac{1}{2}\int_0^T\|y-y_d\|_{L^2(\Omega)}^2dt+\frac{\alpha}{2}\int_0^T\|u\|_{L^2(\Omega)}^2dt
\end{eqnarray}
with the parabolic equation constraint
\begin{eqnarray}
\begin{aligned}\label{eq: constraint}
y_t-\Delta y=&~f+u,~~\mbox{in}~Q_T=\Omega\times(0,T],\\
\mathcal {K}y=&~0,~~\mbox{on}~\partial \Omega,\\
y=&~y_0,~~\mbox{at}~t=0,
\end{aligned}
\end{eqnarray}
where $y_d$ is the desired state and $\alpha>0$ is a
given regularization parameter that also be called the
proportionality factor. $y_0$ is the initial condition, $f$ is the
given source function, and $\mathcal{K}$ is an operator such that
\begin{eqnarray}
\begin{aligned}\label{eq: boundary}
\mathcal {K}y=&~y,~~~~~\mbox{Dirichlet boundary condition},\\
\mathcal {K}y=&~\frac{\partial y}{\partial n},~~\mbox{Neumann
boundary condition}.
\end{aligned}
\end{eqnarray}
The existence and uniqueness of the solution to the
unconstrained distribution control problem \eqref{eq:
objection}-\eqref{eq: constraint} have already been proved in
\cite{Tr10} under some moderate assumptions. The
major contributions of this paper are an efficient method for
solving this problem and the corresponding convergence analysis.
The proposed strategy also could be extended to problems with
observations on part of the domain $\Omega$, or with
boundary controls.

Similar to the traditional optimal control problem constrained with elliptic
equation, there are two major numerical approaches for solving the
parabolic optimal control problems, optimize-then-discretize and
discretize-then-optimize algorithms. The former
approach mainly starts with the continuous Lagrangian function
(\cite{AF12,LLX17})
\begin{eqnarray*}\label{eq: cLagrange}
\mathcal{L}(y,u,p)=\mathcal{J}(y,u)+\int_0^T\big((y_t-f-u,p)+(\nabla
y, \nabla p)\big)dt+\big(y(\cdot,0)-y_0,p(\cdot,0)\big),
\end{eqnarray*}
which implies the first order optimality condition
\begin{eqnarray}
\begin{aligned}\label{eq: optcondition}
y_t-\Delta y=&~f+u,&\mbox{in}~\Omega\times(0,T],\quad \mathcal
{K}y=&~0,&\mbox{on}~\partial \Omega,\quad y=&~y_0,&\mbox{at}~t=0,\\
p_t+\Delta p=&~y-y_d,&\mbox{in}~\Omega\times[0,T),\quad\mathcal
{K}p=&~0,&\mbox{on}~\partial \Omega,\quad p=&~0,&\mbox{at}~t=T,\\
\alpha u-p=&~0, &\mbox{in}~\Omega\times [0,T], \quad~\quad\quad &&
\end{aligned}
\end{eqnarray}
where $p$ is the dual state variable (also called Lagrangian
multiplier). Then discretize the first order optimality condition
\eqref{eq: optcondition} to get the approximations of the state
variable $y$ and the control variable $u$. The existing
discretization methods in spatial direction for dealing with the
optimality condition \eqref{eq: optcondition} include such methods
as the conforming finite element method
(\cite{CGR14,GY16,MV08a,MV08b}), the nonconforming finite
element method (\cite{GS17}), and the finite volume method
\cite{LCHH14}. For the temporal discretization, the key tools are
the finite difference method (\cite{GS17,LCHH14}) and the
discontinuous Galerkin method (\cite{GY16,MV08a,MV08b}). The
best error estimates for these methods are of order $O(h^2+\tau^2)$,
where $h$ and $\tau$ stand for the step size in spatial direction
and temporal direction, respectively. Although the convergence
analysis for the direct discretization methods above has been
established, they are not efficient enough to be practical. The
reason is that the first order optimality condition \eqref{eq:
optcondition} is a forward-backward coupled system, which leads to a
large-scale system and has to be solved simultaneously. Especially,
when the high accuracy or a long time period $T$ is required, the
discretized system is too large to be solved directly. To surmount
this computational challenge, the researchers have developed
techniques such as the multigrid method
(\cite{ADHV11,BW09,LLX17}) and the adaptive method
(\cite{CHY08,HCH11,LMTY04}) for solving the
parabolic optimal control problems. Although these methods perform
efficiently under some moderate assumptions, the design of the multigrid
and the estimation of the posteriori error are still not easy jobs
(\cite{ADHV11,LLX17,LMTY04}).

On the other hand, the latter approach discretize the original
parabolic optimal control problem \eqref{eq: objection}-\eqref{eq:
constraint} directly, and then solve the first order optimality
condition of the discretized optimization problem to obtain the
approximations of the state and control variables. This approach
also needs to overcome the bottle neck issue from the large-scale
discretized system. Based on the traditional
preconditioning techniques for saddle-point systems
that arise from static PDE constrained optimal control problems
(\cite{PW12,PSW12,SU14,YHC16}), some
preconditioning methods are proposed for the time evolution
problems, we refer to \cite{BS15,DSSS13,Kwo17,SB15} and references
therein for the rich literature. Moreover, there also
exist some parallel algorithms for solving the optimal control
problem \eqref{eq: objection}-\eqref{eq: constraint} based on above
two approaches, such as the time domain decomposition methods
\cite{LW18,Ria16}, and the non-intrusive parallel-in-time approach
\cite{GGS19}. The efficiency of these parallel methods have been
proved, but the implements are complicated. The purpose of this
paper is to design a simple and feasible parallel algorithm.

The optimize-then-discretize and discretize-then-optimize approaches
for solving the parabolic optimal control
problems \eqref{eq: objection}-\eqref{eq: constraint} are not
consistent. As a result, the two different approaches could lead to
two different first order optimality discretized systems. The
advantages and disadvantages of both approaches have been summarized
systematically by Apel et. al. (\cite{AF12}) and Gunzburger
(\cite{Gun87}). With the former approach, one can
choose a good approximation of the adjoint equation but the solution
operator may not be symmetric and positive definite. On the other
hand, the discretized gradient is the right direction of descent,
but it is not clear whether the discretized adjoint equation is an
appropriate discretization of the continuous adjoint equation.
Therefore, a consistent scheme which combines the advantages of both
approaches is needed sorely. In \cite{AF12}, Apel et. al.
proposed a consistent Crank-Nicolson finite element method (CN-FEM),
whose convergence rate are of order two in both spatial and temporal
directions, for optimal control problems with evolution equation
constraints. Following this idea, we adopt the CN-FEM to discretize
the parabolic optimal control problem \eqref{eq:
objection}-\eqref{eq: constraint}.

The proposed algorithm in this paper is mainly based on the
discretize-then-optimize approach, which is equivalent to the
optimize-then-discretize approach because of the consistency.
The design of this algorithm is simpler than the existing
methods based on the first order optimality conditions for the original
problem or its discretization form.
In fact, we solve the original discretized optimization problem directly,
rather than the first order optimality condition corresponding to the
discretized optimization problem. After discretization with CN-FEM,
the original problem \eqref{eq: objection}-\eqref{eq: constraint}
can be rewritten as a discretized optimization problem with
separable unknown vectors. Based on the separable structures of the
objective functional and the constraint equations, a full Jacobian
decomposition method with correction (\cite{HHY15}) is proposed to solve the
discretized subsystems in parallel, which improves the computational
efficiency significantly. Moreover, the global convergence analysis
is established, which includes the discretization error by the
CN-FEM and the iteration error by the full Jacobian decomposition
method with correction.

The rest of this paper is organized as follows. In section
\ref{discretization}, the parabolic optimal control problem
\eqref{eq: objection}-\eqref{eq: constraint} is discretized by using
the finite element approximation, the Crank-Nicolson scheme and
the numerical integration formula, and we present
the error estimates between the continuous solutions and their
discretized forms. The full Jacobian decomposition method with correction
for the discretized optimization system, and the corresponding global
error analysis are introduced in section \ref{ajsa}. In section \ref{IMA},
we describe the implementation of the parallel algorithm in details. In
section \ref{NE}, numerical simulations are presented to test the
performance of the proposed method. The last section is devoted to
some concluding remarks.



\section{Finite element method and Crank-Nicolson scheme}\label{discretization}

In this section, the linear finite element is applied to
approximate the state variable and the control variable
in spatial direction. Moreover, in temporal direction,
the constraint equation \eqref{eq: constraint} is discretized
by the Crank-Nicolson scheme, and the integrals with respect
to the state variable and the control variable in the objective
function \eqref{eq: objection} are approximated by the trapezoidal
rule and midpoint rule, respectively. Now, some notations are introduced
which shall be used in the sequel. Define the spaces
\begin{eqnarray*}
H^{1}(\Omega)&:=&\{v\in L^2(\Omega)~\big|~ \frac{\partial
v}{\partial
x_j}\in L^2(\Omega), j=1,2\},\\
H^{1}_{0}(\Omega)&:=&\{v\in H^1(\Omega)~\big|~ v=0~on~
\partial \Omega\}.
\end{eqnarray*}
We denote $V=H_0^1(\Omega)$ and $V=H^1(\Omega)$ for Dirichlet
boundary condition and Neumann boundary condition, respectively,
whenever there is no ambiguity. Denote by $V^*$ the dual space of
$V$. Let $X$ be a Banach space with norm $\|\cdot\|_X$. For each
$T>0$, we denote by $L^2(0,T;X)$ the $X$-valued $L^2$ space
consisting of strongly measurable functions $v : (0, T] \rightarrow
X$ such that
\begin{eqnarray*}
\|v\|_{L^2(0,T;X)}:= \left(\int_0^T
\|v(\cdot,t)\|^2_{X}dt\right)^{\frac{1}{2}} <\infty.
\end{eqnarray*}

Firstly, we disrectize the parabolic equation constraint based on
the variational formulation. For any $u\in L^2(0,T;V)$, the
variational formulation of \eqref{eq: constraint} is given by: Find
$y\in L^2(0,T;V)$, such that $y(x,0)=y_0$ and
\begin{eqnarray}\label{eq:CVI}
(\mbox{\textbf{VF}})\qquad (y_t, v)+(\nabla y, \nabla
v)=(f+u,v),~~\forall~ v\in V, 0 < t \leq T.
\end{eqnarray}

Secondly, we present the semi-discrete approximation of the
variational formulation (\ref{eq:CVI}). Let $\mathcal{J}_t: 0 = t_0
< t_1 < \dots < t_M = T$ be an equidistant partition on $[0, T]$
with $\tau=\frac{T}{M}$ standing for the temporal step size. Assume
that $\mathcal{T}_h$ is a shape regular triangulation on $\Omega$
consisting of polygons. For any $K\in \mathcal{T}_h$, $h_K$ stands
for the diameter of the polygon $K$ and $h = \max\limits_{K\in
\mathcal{T}_h} h_K$. Assume that $P_{T_b}$ and $P_{T_i}$ are the
sets of the grid points on the boundary and in the interior,
respectively. Let $N_b = \sharp\{P_{T_b}\}$, $N_i =
\sharp\{P_{T_i}\}$, and $N = N_b + N_i$. Define the piecewise linear
element space
\begin{eqnarray*}\label{eq:space Sh}
 S^1(\Omega):&=&\{v~\big|~ v\in  H^1(\Omega),~v|_{K} \in \mathcal{P}_{1}, \forall~K \in \mathcal{T}_h\},\\
 S^1_0(\Omega):&=&\{v~\big|~ v\in H^1_0(\Omega),~v|_{K} \in \mathcal{P}_{1}, \forall~K \in \mathcal{T}_h\},
\end{eqnarray*}
where $\mathcal{P}_{1}$ stands for the set of polynomials with
degree less than or equal to one. Similarly, we denote
$V_h=S^1_0(\Omega)$ and $V_h=S^1(\Omega)$ for Dirichlet boundary
condition and Neumann boundary condition, respectively, whenever
there is no ambiguity. Then for any $u_h \in L^2(0,T;V_h)$, the
semi-discrete approximation of the variational formulation
(\ref{eq:CVI}) is: Find $y_h\in L^2(0,T;V_h)$, such that
$y_h(x,0)=Q_hy_0$ and
\begin{eqnarray}\label{eq:CVI-SD}
(y_{ht}, v_h)+(\nabla y_h, \nabla v_h)=(u_h,v_h),~~\forall~ v_h\in
V_h, 0<t\leq T,
\end{eqnarray}
where $Q_h$ is the finite element interpolation operator.

Let $S^1_0(\Omega)=\mbox{span}\{\varphi_{i_j},j=1,\cdots,N_i\}$ and
$S^1(\Omega)=\mbox{span}\{\varphi_{i_j},j=1,\cdots,N_i;$ $\varphi_{b_j},j=1,\cdots,N_b\}$,
where $\varphi_{i_j}$ and $\varphi_{b_j}$ are the basis functions
corresponding to the grid points $g_{i_j}\in P_{T_i}$ and
$g_{b_j}\in P_{T_b}$, respectively. At each point $t=t_m
(m=1,\dots,M)$, the finite element approximation of the functions
$y(x,t_m)$ and $u(x,t_m)$ are given by
\begin{eqnarray*}
\begin{aligned}
 y^m_h&=\sum\limits_{j=1}^{N_i}y^m_{i_j}\varphi_{i_j}(x), &&u^m_h=\sum\limits_{j=1}^{N_i}u^m_{i_j}\varphi_{i_j}(x),
 &&\mbox{Dirichlet case},\\
 y^m_h&=\sum\limits_{j=1}^{N_i}y^m_{i_j}\varphi_{i_j}(x)+\sum\limits_{j=1}^{N_b}y^m_{b_j}\varphi_{b_j}(x),
 &&u^m_h=\sum\limits_{j=1}^{N_i}u^m_{i_j}\varphi_{i_j}(x)+\sum\limits_{j=1}^{N_b}u^m_{b_j}\varphi_{b_j}(x),
 &&\mbox{Neumann case}.
\end{aligned}
\end{eqnarray*}
Thirdly, the fully discretized approximation based on the Crank-Nicolson
scheme is given by
\begin{eqnarray}\label{eq:CVI-FD}\hspace{-5mm}
(\frac{y^{m+1}_{h}-y^{m}_{h}}{\tau}, v_h)+(\frac{\nabla
y^{m+1}_h+\nabla y^{m}_h}{2}, \nabla v_h)=(f^{m+\frac{1}{2}}+
u^{m+\frac{1}{2}}_h,v_h),~\forall~ v_h\in V_h, m=1,\dots,M,
\end{eqnarray}
where $f^{m+\frac{1}{2}}=f(x,\frac{t_m+t_{m-1}}{2})$ and
$u^{m+\frac{1}{2}}_h$ is the finite element approximation of the
function $u(x,\frac{t_m+t_{m-1}}{2})$.
The corresponding matrix-vector form of the Crank-Nicolson scheme
\eqref{eq:CVI-FD} is
\begin{eqnarray}\label{eq:CVI-M}
(A+\frac{\tau}{2}B)Y_{m+1}-(A-\frac{\tau}{2}B)Y_{m}-\tau A
U_{m+\frac{1}{2}}=F_{m+\frac{1}{2}}, \qquad m=1,\cdots,M,
\end{eqnarray}
where
\begin{eqnarray*}
\begin{aligned}
&A_{k,j}=(\varphi_{i_j}(x),\varphi_{i_k}(x)),~~B_{k,j}=(\nabla\varphi_{i_j}(x),\nabla\varphi_{i_k}(x)),~~j,k=1,2\cdots,N_i,\\
&Y_m=(y^m_{i_1},y^m_{i_2},\cdots,y^m_{i_{N_i}})^T,~~U_{m+\frac{1}{2}}=(u^{m+\frac{1}{2}}_{i_1},u^{m+\frac{1}{2}}_{i_2},\cdots,
u^{m+\frac{1}{2}}_{i_{N_i}})^T,\\
&(F_{m+\frac{1}{2}})_k=\tau(f^{m+\frac{1}{2}},\varphi_{i_k}(x)),~~k=1,2\cdots,N_i,
\end{aligned}
\end{eqnarray*}
in case of Dirichlet boundary condition, and
\begin{eqnarray*}\hspace{-5mm}
\begin{aligned}
&A=\left[\begin{array}{cc}
            A^{ii}   &   A^{bi}    \\[0.4cm]
            A^{ib}   &   A^{bb}
 \end{array}\right],\qquad
 B=\left[\begin{array}{cc}
            B^{ii}   &   B^{bi}    \\[0.4cm]
            B^{ib}   &   B^{bb}
 \end{array}\right], \\
\end{aligned}
\end{eqnarray*}
\begin{eqnarray*}\hspace{-5mm}
\begin{aligned}
&A^{ii}_{k,j}=(\varphi_{i_j}(x),\varphi_{i_k}(x)),~~B^{ii}_{k,j}=(\nabla\varphi_{i_j}(x),\nabla\varphi_{i_k}(x)),~~k=1,2\cdots,N_i, ~~j=1,2\cdots,N_i,\\
&A^{ib}_{k,j}=(\varphi_{i_j}(x),\varphi_{b_k}(x)),~~B^{ib}_{k,j}=(\nabla\varphi_{i_j}(x),\nabla\varphi_{b_k}(x)),~~k=1,2\cdots,N_b, ~~j=1,2\cdots,N_i,\\
&A^{bi}_{k,j}=(\varphi_{b_j}(x),\varphi_{i_k}(x)),~~B^{bi}_{k,j}=(\nabla\varphi_{i_j}(x),\nabla\varphi_{i_k}(x)),~~k=1,2\cdots,N_i, ~~j=1,2\cdots,N_b,\\
&A^{bb}_{k,j}=(\varphi_{b_j}(x),\varphi_{b_k}(x)),~~B^{bb}_{k,j}=(\nabla\varphi_{i_j}(x),\nabla\varphi_{b_k}(x)),~~k=1,2\cdots,N_b,
~~j=1,2\cdots,N_b,\\
&Y_m=(y^m_{i_1},\cdots,y^m_{i_{N_i}},y^{m}_{b_1},\cdots,y^{m}_{b_{N_b}})^T,\quad
U_{m+\frac{1}{2}}=(u^{m+\frac{1}{2}}_{i_1},\cdots,
u^{m+\frac{1}{2}}_{i_{N_i}},u^{m+\frac{1}{2}}_{b_1},\cdots,u^{m+\frac{1}{2}}_{b_{N_b}})^T,\\
&(F_{m+\frac{1}{2}})_k=\tau(f^{m+\frac{1}{2}},\varphi_{i_k}(x)),
k=1,2\cdots,N_i,~~
(F_{m+\frac{1}{2}})_{N_i+j}=\tau(f^{m+\frac{1}{2}},\varphi_{b_j}(x)),~j=1,2\cdots,N_b,
\end{aligned}
\end{eqnarray*}
in case of Neumann boundary condition. Furthermore, the Crank-Nicolson
discretization of the parabolic equation \eqref{eq: constraint} can
be rewritten as
\begin{eqnarray}\label{eq:CVI-globleM}
\mathcal{A}Y+\mathcal{B}U=\mathcal{F},
\end{eqnarray}
where
\begin{eqnarray*}
&&\mathcal{A}=\left[
\begin{array}{ccccc}
        A+\frac{\tau}{2}B&0&0&\cdots&0 \\
        -A+\frac{\tau}{2}B&A+\frac{\tau}{2}B&0&\cdots&0 \\
        0&-A+\frac{\tau}{2}B&\ddots&\ddots& \vdots \\
        \vdots&\ddots&\ddots&A+\frac{\tau}{2}B&  0 \\
        0&\cdots&0&-A+\frac{\tau}{2}B& A+\frac{\tau}{2}B
        \end{array}\right],
\quad Y=\left[
\begin{array}{c}
        Y_1 \\
        Y_2 \\
        \vdots \\
        Y_{M-1} \\
        Y_{M}
        \end{array}\right],
\end{eqnarray*}
\begin{eqnarray*}
&&\mathcal{B}=\left[
\begin{array}{ccccc}
        -\tau A&0&0&\cdots&0 \\
        0&-\tau A&0&\cdots&0 \\
        0&0&\ddots&\ddots& \vdots \\
        \vdots&\ddots&\ddots&-\tau A&  0 \\
        0&\cdots&0&0&-\tau A
        \end{array}\right],
\quad U=\left[
\begin{array}{c}
        U_{\frac{1}{2}} \\
        U_{\frac{3}{2}} \\
        \vdots \\
        U_{M-\frac{3}{2}} \\
        U_{M-\frac{1}{2}}
        \end{array}\right],
\quad \mathcal{F}=\left[
\begin{array}{c}
        F_{\frac{1}{2}}+(A-\frac{\tau}{2}B)Y_0 \\
        F_{\frac{3}{2}}\\
        F_{\frac{5}{2}}\\
        \vdots \\
        F_{M-\frac{1}{2}}
        \end{array}\right].
\end{eqnarray*}
Let $\mathcal{A}_m$ and $\mathcal{B}_m$ be the block columns,
corresponding to $t_m$, of $\mathcal{A}$ and $\mathcal{B}$,
respectively, then the Crank-Nicolson
discretization \eqref{eq:CVI-globleM} is equivalent to the following
formulation with separable structures
\begin{eqnarray}\label{eq:CVI-SglobleM}
\sum_{m=1}^M \mathcal{A}_m Y_m+\sum_{m=1}^M
\mathcal{B}_mU_{m-\frac{1}{2}}=\mathcal{F}.
\end{eqnarray}

Next, we discretize the objective functional \eqref{eq: objection}.
In the temporal direction, the integrations with respect to the
state variable $y$ and control variable $u$ are discretized by the
trapezoidal rule and midpoint rule, respectively, which are
consistent with the Crank-Nicolson scheme \eqref{eq:CVI-FD}.
Meanwhile, we still use finite element approximation in the spatial
direction. This discretion approach guarantees the consistency of
the discretization and the optimization \cite{ADHV11}. Therefore,
the discretization formulation of objective functional \eqref{eq:
objection} is given by
\begin{eqnarray}\label{eq: FEMobjection}
\begin{aligned}
\mathcal{J}(y_h,u_h)=&\,\frac{1}{2}\left(\frac{\tau}{2}\|y_h^0-y_{d}^0\|_{L^2(\Omega)}^2+\tau
\sum_{m=1}^{M-1}\|y_h^m-y_{d}^m\|_{L^2(\Omega)}^2+\frac{\tau}{2}\|y_h^M-y_{d}^M\|_{L^2(\Omega)}^2\right)\\
&\,+\frac{\alpha \tau}{2}
\sum_{m=0}^{M-1}\|u_h^{m+\frac{1}{2}}\|_{L^2(\Omega)}^2,
\end{aligned}
\end{eqnarray}
which can be rewritten as the following vector form
\begin{eqnarray}\label{eq: MVobjection}
\begin{aligned}
\hspace{-0.5cm}\mathcal{J}_h(Y,U)=\frac{\tau}{2}\sum\limits_{m=1}^{M-1}\left(Y_m^TAY_m-2d_m^TY^m\right)+\frac{\tau}{4}\left(Y_M^TAY_M-2d_M^TY_M\right)
+\frac{\alpha\tau}{2}\sum\limits_{m=0}^{M-1}U_{m+\frac{1}{2}}^TAU_{m+\frac{1}{2}},
\end{aligned}
\end{eqnarray}
where $d_{j}^m=(\varphi_{i_j}(x),y_d(\cdot,t_m))$, $j=1,\dots,N_i$
in case of Dirichlet boundary condition, and
$d_{j}^m=(\varphi_{i_j}(x),$ $y_d(\cdot,t_m)), j=1,\dots,N_i, d_{N_i+j}^m=(\varphi_{b_j}(x),y_d(\cdot,t_m)), j=1,\dots,N_b$ in
case of Neumann boundary condition.

Finally, the optimal control problem \eqref{eq: objection}-\eqref{eq:
constraint} can be approximated by the discretized system \eqref{eq:
FEMobjection} and \eqref{eq:CVI-FD}, or their vector forms
\eqref{eq: MVobjection} and \eqref{eq:CVI-SglobleM},
which is a large-scale quadratic optimization problem
with linear constraint.

Let
$\|\cdot\|_{L^2(Q_T)}=\left(\displaystyle\int_0^T\|\cdot\|_{L^2(\Omega)}^2dt\right)^{\frac{1}{2}}$
and $\mathbf{\Pi_t}$ be the piecewise linear interpolation operator
in temporal direction. Under some moderate assumptions, Apel and Flaig
have obtained the following convergence estimates for discretized
system \eqref{eq: FEMobjection} and \eqref{eq:CVI-FD}.
\begin{lem}(cf.~\cite{AF12})\label{lem:fem}
Let $(y^{*},~p^{*},~u^{*})$ be the exact solution of the first order
optimality condition \eqref{eq: optcondition}, and
$\{y_h^m,~u_h^{m-\frac{1}{2}}\}_{m=1}^{M}$ be the optimal solution
of the discretized optimization problem \eqref{eq: FEMobjection} and
\eqref{eq:CVI-FD}. Assume $y^{*}$, $p^{*}$, $u^{*}$, $f\in
H^3(0,T;L^2(\Omega))\cap H^2(0,T;H^2(\Omega))$ and $y_0\in
H^2(\Omega)$, then
\begin{eqnarray}\label{eq:fem}
\begin{aligned}
\|y_h^m-y^{*}(\cdot,t_m)\|_{L^2(\Omega)}\leq &~C_1 h^2+C_2 \tau^2, \qquad m=1,\cdots,M,\\
\|\mathbf{\Pi_t}u_h-u^{*}\|_{L^2(Q_T)}\leq &~C_1 h^2+C_2 \tau^2,
\end{aligned}
\end{eqnarray}
where $C_1$ and $C_2$ are constants independent of $h$ and $\tau$.
\end{lem}

\section{Full Jacobian decomposition algorithm with correction }\label{ajsa}

In this section, the full Jacobian decomposition method with correction
is applied to solve the optimization problem \eqref{eq: MVobjection}
with the constraint \eqref{eq:CVI-SglobleM}. Different from the
existing discretize-then-optimize methods that are used to deal with
the first order optimality condition of the optimization problem
(cf. \cite{PSW12}), we solve the original optimization problem by
the full Jacobian iteration with correction directly, which could avoid
solving the large-scale coupled system, by parallel computing.

For the convenience of the expression, we first reformulate the
optimization problem \eqref{eq: MVobjection} with the linear
constraint \eqref{eq:CVI-SglobleM} as follows

\begin{eqnarray}\label{eq: A-Problem-M}
 \begin{array}{rl}
     \min\limits_{z_l\in \mathcal{Z}_{l}} & \sum\limits_{l=1}^{2M} \theta_l(z_l)    \\[0.2cm]
    s.t.~ & \sum\limits_{l=1}^{2M} \mathcal{M}_l z_l =\mathcal{F},
          \end{array}
\end{eqnarray}
where $\mathcal{Z}_{l}=\mathbb{R}^{N_i}$ for Dirichlet case,
$\mathcal{Z}_{l}=\mathbb{R}^{N}$ for Neumann case, and

\begin{eqnarray*}
\begin{aligned}
z_{2l-1}=&~U_{l-\frac{1}{2}},~~z_{2l}=Y_l,~~l=1,\cdots, M,\\
\theta_{2l-1}(z_{2l-1})=&~\frac{\alpha\tau}{2}U_{l-\frac{1}{2}}^TAU_{l-\frac{1}{2}},~~l=1,\cdots, M,\\
\theta_{2l}(z_{2l})=&~\frac{\tau}{2}(Y_l^TAY_l-2d_l^TY_l),~~l=1,\cdots, M-1,\\
\theta_{2M}(z_{2M})=&~\frac{\tau}{4}(Y_M^TAY_M-2d_M^TY_M),\\
\mathcal{M}_{2l-1}=&~\mathcal{B}_{l},~~\mathcal{M}_{2l}=\mathcal{A}_{l},~~l=1,\cdots,M.
\end{aligned}
\end{eqnarray*}

The augmented Lagrangian method
(ALM) is an efficient and robust algorithm for solving the above optimization problem (cf. \cite{HHY15}).
Let the Lagrangian function of the optimization problem \eqref{eq: A-Problem-M} be
\begin{equation}\label{eq:Lag}
     L(z_1, \ldots, z_{2M}, \lambda) = \sum_{l=1}^{2M} \theta_l(z_l)
     -\lambda^T\bigl(\sum_{l=1}^{2M} \mathcal{M}_l z_l -\mathcal{F}\bigr),
   \end{equation}
and the corresponding augmented Lagrangian function be
\begin{equation}\label{eq:Alag}
L_A (z_1, \ldots, z_{2M}, \lambda) =  L(z_1, \ldots, z_{2M},
\lambda)+  \frac{\beta}{2}\Big\|\sum_{l=1}^{2M} \mathcal{M}_l z_l
-\mathcal{F}\Big\|^2,
\end{equation}
where the Lagrange multiplier $\lambda \in \mathbb{R}^{MN_i}$ for
Dirichilet case (or $\mathbb{R}^{MN}$ for Neumann case) and the
positive penalty parameter $\beta \in \mathbb{R}$. Here and
hereafter, $\|\cdot\|$ stands for the $l_2$ norm of the vectors.
Applying the ALM scheme directly to the well-structured form
\eqref{eq: A-Problem-M}, we obtain the following the iterative
scheme
\begin{eqnarray}\label{eq:DALM}
\left\{\begin{array}{rcl} (z_1^{k+1}, \ldots, z_{2M}^{k+1}) &=&
\arg\min\bigl\{L_A (z_1, \ldots, z_{2M}, \lambda^k) \, \big|\,  z_l\in \mathcal{ Z}_l,\; l=1,\cdots, 2M\bigr\},\\[0.2cm]
 \lambda^{k+1} &=&\lambda^k-\beta(\sum\limits_{l=1}^{2M} \mathcal{M}_l z_l^{k+1} -\mathcal{F}).
\end{array}\right.
\end{eqnarray}
The convergence of the sequence generated by \eqref{eq:DALM} is
well-known. But there is a huge challenge in
implementation when the number of the vectors is greater than $2$.
It is because that all the subvectors $z_l$ are required to be
solved simultaneously and all $\theta_l$ have to be considered
aggregately. Taking advantage of the separable structure of
objective function and constraints, we decouple the variable $z$
into $2M$ components. Meanwhile the objective function is decomposed
into $2M$ components, where the $l$-th component only involves
$\theta_l(z_l)$ and certain quadratic polynomials in $z_l$, which
leads to $2M$ subproblems. Furthermore, one can easily deduce the
closed-form solution for each subproblem. This kind of splitting
techniques are widely used in many applications arising from diverse
areas such as image processing, statistical learning, and
compressive sensing. One of the important versions is the ALM with
full Jacobian decomposition, the corresponding subproblems for
solving (\ref{eq:DALM}) by this splitting method are given by:

\begin{eqnarray}\label{ADMM-GP}\hspace{8mm}
\left\{\!
\begin{array}{l}
z_1^{k+1} =\arg\min
 \bigl\{\theta_1(z_1)-z_1^T\mathcal{M}_1^T\lambda^k + { \textstyle{\frac{\beta}{2}}}\|\mathcal{M}_1z_1 + \sum\limits_{j=2}^{2M}
\mathcal{M}_jz_j^k-\mathcal{F} \|^2\; \big| \; z_1\in\mathcal{Z}_1  \bigr\};\\[0.2cm]
z_2^{k+1} =\arg\min
\bigl\{ \theta_2(z_2)-z_2^T\mathcal{M}_2^T\lambda^k + {\textstyle{\frac{\beta}{2}}}\|\mathcal{M}_1z^k_1 + \mathcal{M}_2z_2 +\sum\limits_{j=3}^{2M} \mathcal{M}_jz_j^k-\mathcal{F} \|^2 \;\big| \; z_2\in\mathcal{Z}_2  \bigr\};\\
     \cdots \cdots \\
z_l^{k+1} =\arg\min \bigl\{\theta_l(z_l)\!-\!
z_l^T\mathcal{M}_l^T\lambda^k\! +\!
{\textstyle{\frac{\beta}{2}}}\|\sum\limits_{j=1}^{l-1}\!
\mathcal{M}_j z_j^k
\!+ \!\mathcal{M}_lz_l \! + \!\sum\limits_{j=l+1}^{2M}\! \mathcal{M}_jz_j^k-\mathcal{F}\|^2\;\big|\; z_l \in\mathcal{Z}_l \bigr\};\\
      \cdots \cdots \\
z_{2M}^{k+1}=\!\arg\min \bigl\{ \theta_{2M}(z_{2M})
-z_{2M}^T\mathcal{M}_{2M}^T\lambda^k +
{\textstyle{\frac{\beta}{2}}}\|\sum\limits_{j=1}^{2M} \mathcal{M}_j
z_j^k + \mathcal{M}_{2M}z_{2M}-\mathcal{F}\|^2
\;\big| \; z_{2M} \in\mathcal{Z}_{2M} \bigr\};\\[0.5cm]
{\lambda}^{k+1} = \lambda^k -\beta(\sum\limits_{j=1}^{2M}\mathcal{M}_j
z_j^{k+1}-\mathcal{F}).
\end{array}
\right.
\end{eqnarray}

The splitting version of ALM with full Jacobian decomposition
(\ref{ADMM-GP}) allows all the $z_l$-subproblems being solved in
parallel, and this is an extremely important feature when large or
huge scale data are considered with parallel computing
infrastructures being available. But the Jacobian splitting scheme
(\ref{ADMM-GP}) is not convergent \cite{CHYY16}. Fortunately, based
on it, Yuan et. al. \cite{HHY15} proposed a full Jacobian
decomposition method with correction, which is convergent,
and the similar idea also be introduced in \cite{BK19}.
They use the output of (\ref{ADMM-GP}), denoted
by ${\widetilde w}^k:=({\widetilde z}_1^k, {\widetilde z}_2^k,
\cdots,{\widetilde z}_{2M}^k, {\widetilde \lambda^k})$, as a
predictor, and added a correct step to update the iteration solution
$w^{k+1}$. We use the same notations in our paper as well.
The splitting version of ALM with full
Jacobian decomposition and a corrector with constant step size can
be described as\vspace{5mm}

\medskip
\noindent \fbox{\begin{minipage}{12.5cm}
  \noindent {{\bf Splitting Algorithm 1}}

\indent Step 1: Generate $\widetilde{w}^k$ via (\ref{ADMM-GP}).

\indent Step 2: Generate the new iterate
$w^{k+1}$ via

\begin{subequations} \label{Cor-2}
\[ \label{Cor-2w}
   w^{k+1}  = w^k- \nu(w^k-\widetilde{w}^k),
      \]
where
\[ \label{Cor-2a}
     \nu= \gamma\bigl(1-\sqrt{\textstyle\frac{2M}{2M+1}}\,\bigr) \qquad \hbox{and}
      \qquad \gamma\in(0,2).   \]
\end{subequations}

\end{minipage}}
\vspace{7mm}

Yuan et. al. showed the contraction property of the above splitting
algorithm, and obtained the convergence of the iterations under some
moderate assumptions \cite{HHY15}.
\vspace{5mm}

\begin{lem}\label{lem: contraction}
Let $\{w^k\}$ be the sequence generated by the splitting algorithm 1
with an arbitrary initial iterate $w^0$, and $w^{*}$ be the saddle
point of the Lagrange function  \eqref{eq:Lag}. Then,
\begin{eqnarray}\label{eq: contraction}
\begin{aligned}
\|w^{k+1}-w^{*}\|_H^2\leq&~\|w^k-w^{*}\|_H^2-\frac{2-\gamma}{\gamma}\|w^k-w^{k+1}\|_H^2,\\
\|w^k-w^{k+1}\|_H^2\leq&~\frac{4}{\gamma(2-\gamma)(k+1)}\|w^{0}-w^{*}\|_H^2,
\end{aligned}
\end{eqnarray}
where

\[   \label{Matrix-H}
    H =\beta
   \left(\begin{array}{ccccc}
            2\mathcal{M}_1^T\mathcal{M}_1   &   \mathcal{M}_1^T\mathcal{M}_2    &   \cdots      & \mathcal{M}_1^T\mathcal{M}_{2M}             &   0    \\[0.4cm]
            \mathcal{M}_2^T\mathcal{M}_1    &      \ddots                       &   \ddots      &  \vdots                                     & \vdots \\[0.4cm]
            \vdots                          &      \ddots                       &   \ddots      &  \mathcal{M}_{2M}^T\mathcal{M}_{2M} &  \vdots \\[0.4cm]
            \mathcal{M}_{2M}^T\mathcal{M}_1 & \cdots                          &\!\! \mathcal{M}_{2M}^T\mathcal{M}_{2M} &  2\mathcal{M}_{2M}^T\mathcal{M}_{2M}&   0   \\[0.4cm]
           0                                & \cdots                            &  \cdots               &    0               &  \dfrac{1}{\beta^2}
 \end{array}\right).
 \]
\end{lem}

\begin{lem}\label{lem: convergence}
Under the conditions in Lemma \ref{lem: contraction}, if the
columns $\mathcal{M}_l\,(l =1,\dots,2M)$ in \eqref{eq: A-Problem-M}
are all of full column rank, then $\{w^k\}$ converges to the saddle
point $w^{*}$ of the Lagrange function \eqref{eq:Lag}.
\end{lem}

\begin{thm}\label{thm: ALMconverge}
The sequence $\{w^k\}$ generated by the splitting algorithm 1
converges to the saddle point $w^{*}$ of the Lagrange function
\eqref{eq:Lag}.
\end{thm}
\textbf{Proof.} By Lemma \ref{lem: convergence}, we only need to
prove that $\mathcal{M}_l$ ($l =1,\cdots,2M$) in \eqref{eq:
A-Problem-M} are all of full column rank, which is equivalent to
proving that $\mathcal{A}_m$ and $\mathcal{B}_m$ ($m =1,\cdots,M$)
in \eqref{eq:CVI-SglobleM} are of full column rank.

From the approximation \eqref{eq:CVI-M}, we know that the mass
matrix $A$ and stiffness matrix $B$ formed by the linear elements
are all positive definite, which means $A$ and $A+\frac{\tau}{2}B$
have full column rank. This, together with the definitions of
$\mathcal{A}$ and $\mathcal{B}$ in \eqref{eq:CVI-globleM} implies
the conclusion.   $\hfill{} \Box$

Let $\mathbf{R_x}$ and $\mathbf{\Pi_t}$ be the linear interpolation
operators in spatial direction and temporal direction, respectively.
Based on Lemma \ref{lem:fem}, Theorem \ref{thm: ALMconverge}, and
the fact that $w^{*}$ is the saddle point of the Lagrange function
\eqref{eq:Lag} if and only if
$(U_{\frac{1}{2}}^*,Y_{1}^*\cdots,U_{M-\frac{1}{2}}^*,Y_{M}^*)$ is
the optimal solution of the optimization problem \eqref{eq:
A-Problem-M}, we obtain the following convergence result.

\begin{thm}\label{thm: Allconverge}
Let $\{u^*$, $y^*\}$ be the optimal solution of the optimal control problem
\eqref{eq: objection}-\eqref{eq: constraint} and
$w^k=(U_{\frac{1}{2}}^k,Y_{1}^k\cdots,U_{M-\frac{1}{2}}^k,Y_{M}^k,\lambda^k)$
be the sequence generated by the splitting algorithm 1. Then,
$\mathbf{R_x} Y_m^{k}$ converges to the optimal state $y^*(x,t_m)$,
and $\mathbf{\Pi_t} \mathbf{R_x} \{U_{m-\frac{1}{2}}^k\}$ converges
to the optimal control $u^*(x,t)$ as $(h, \tau,
\frac{1}{k})\rightarrow(0,0, 0)$.
\end{thm}
\textbf{Proof.} Using the identity $y_h^m= \mathbf{R_x} Y_m^*$ and triangle inequality, we obtain
\begin{eqnarray}\label{eq: yconverge}
\begin{aligned}
\|\mathbf{R_x}Y_{m}^k-y^*(x,t_m)\|_{L^2(\Omega)}\leq&~\|\mathbf{R_x}Y_{m}^k-\mathbf{R_x}Y_{m}^*\|_{L^2(\Omega)}
+\|\mathbf{R_x}Y_{m}^*-y^*(x,t_m)\|_{L^2(\Omega)},\\
\leq&~C_0 h\|Y_{m}^k-Y_{m}^*\|+\|y_h^m-y^*(x,t_m)\|_{L^2(\Omega)},\\
\leq&~C_0 h\|w^{k}-w^{*}\|+C_1 h^2+C_2 \tau^2.\\
\end{aligned}
\end{eqnarray}
Similarly, using the identity $\mathbf{\Pi_t}u_h=
\mathbf{\Pi_t}\mathbf{R_x} \{U_{m-\frac{1}{2}}^*\}$ and triangle
inequality, we obtain
\begin{eqnarray}\label{eq: uconverge}
\begin{aligned}
&~\|\mathbf{\Pi_t}\mathbf{R_x}
\{U_{m-\frac{1}{2}}^k\}-u^*(x,t)\|_{L^2(Q_T)}\\
\leq&~\left\|\mathbf{\Pi_t}\mathbf{R_x}
\{U_{m-\frac{1}{2}}^k\}-\mathbf{\Pi_t}\mathbf{R_x}
\{U_{m-\frac{1}{2}}^*\}\right\|_{L^2(Q_T)}+\|\mathbf{\Pi_t}\mathbf{R_x}
\{U_{m-\frac{1}{2}}^*\}-u^*(x,t)\|_{L^2(Q_T)},\\
\leq&~C_0 h\sqrt{\tau}\|U^{k}-U\|+\|\mathbf{\Pi_t}u_h-u^*(x,t)\|_{L^2(Q_T)},\\
\leq&~C_0 h\sqrt{ \tau}\|w^{k}-w^{*}\|+C_1 h^2+C_2 \tau^2,
\end{aligned}
\end{eqnarray}
where $Q_T=\Omega\times(0,T)$. The overall convergent results follows
from inequality \eqref{eq: yconverge}, \eqref{eq: uconverge}, Lemma
\ref{lem:fem}, Theorem \ref{thm: ALMconverge}, and the equivalence
of norm $\|\cdot\|$ and $\|\cdot\|_H$ directly. $\hfill{} \Box$

\begin{rem}
The estimates \eqref{eq: contraction} is regarded as a testimony of
the theoretical convergence rate $O(\frac{1}{k})$ for the splitting
algorithm 1, which implies $\|w^{k}-w^{*}\|^2_H =O(\frac{1}{k})$ in
general case \cite{HHY15}. The computational results always show
faster convergent rates than $O(\frac{1}{k})$ (see numerical
simulations in section \ref{NE}).
If we set $\tau=h$, based on the definition of $H$,
It has $\|w^{k}-w^{*}\|_H =O(h)\|w^{k}-w^{*}\|$.
Using similar arguments as in the
proof of Theorem \ref{thm: Allconverge}, we can obtain the following
results
\begin{eqnarray}\label{eq: yuconverge}
\begin{aligned}
\|\mathbf{R_x}Y_{m}^k-y^*(x,t_m)\|_{L^2(\Omega)}=&~O(\frac{1}{\sqrt{k}})+O (h^2),\\
\|\mathbf{\Pi_t}\mathbf{R_x}\{U_{m-\frac{1}{2}}^k\}-u^*(x,t)\|_{L^2(Q_T)}
=&~O(\sqrt{\frac{h}{k}})+O (h^2).\\
\end{aligned}
\end{eqnarray}
\end{rem}

\section{Parallel implementation of our algorithm}\label{IMA}

In this section, we present the parallel implementation of
the proposed method in details. Observing the $l$-th
subproblem at the stage of the splitting step (\ref{ADMM-GP}),
we can find that this subproblem is an unconditional
extremum problem only related to $z_l$. Therefore,
the first-order optimality conditions of the subproblems are given by
\begin{eqnarray}\label{eq: firstc}
\theta'_l(\widetilde{z}_l^k) -\mathcal{M}_l^T\lambda^k  +
\beta\mathcal{M}_l^T\bigl(\mathcal{M}_l\widetilde{z}_l^k +
\sum_{j=1,\,j\ne l}^{2M}\mathcal{M} _jz_j^k-\mathcal{F}
\bigr)=0,~~l=1,2,\cdots,2M.
\end{eqnarray}
Now, we consider the iterations on the control variables firstly.
Letting $q^k= \sum\limits_{l=1}^{2M} \mathcal{M}_l
z_l^k-\mathcal{F}-\dfrac{\lambda^k}{\beta}$ yields
\begin{eqnarray*}
(\alpha \tau A+\beta\tau^2 A^TA)\widetilde{z}_l^k +
\beta\mathcal{M}_l^T\bigl(
         q^k-\mathcal{M}_lz_l^k \bigr)=0,~~l=1,3\cdots, 2M-1.
\end{eqnarray*}
Let $q^k=(q_1^k,q_2^k,\dots,q_{M}^k)^T$, we have
\begin{eqnarray*}
\begin{aligned}
&(\alpha \tau A+\beta\tau^2 A^TA)\widetilde{z}_l^k + \beta\tau
A^T\bigl(
         q_{\frac{l+1}{2}}^k-\tau A z_l^k \bigr)=0,~~l=1,3\cdots, 2M-1.
\end{aligned}
\end{eqnarray*}
For simplification, we set $D_{u,z}=\tau^2 A^TA$ and $D_{u,q}=\tau
A$. Then the first-order optimality conditions of the subproblems for
the control variables are given by
\begin{eqnarray}\label{ADMMu}
\begin{aligned}
&(\alpha D_{u,q}+\beta D_{u,z})(\widetilde{z}_1^k,
\widetilde{z}_3^k, \cdots, \widetilde{z}_{2M-1}^k)\\
=&\beta \left(D_{u,z}(z_1^k, z_3^k, \cdots,
z_{2M-1}^k)-D_{u,q}(q_1^k, q_2^k, \cdots, q_{M}^k)\right).
\end{aligned}
\end{eqnarray}

Next, we discuss the state variables. Similar to the control variables, it follows from \eqref{eq: firstc} that
\begin{eqnarray*}
\begin{aligned}
&(\tau A+\beta \mathcal{M}_l^T \mathcal{M}_l)\widetilde{z}_l^k +
\beta\mathcal{M}_l^T\bigl(
         q^k-\mathcal{M}_lz_l^k \bigr)-\tau d_{\frac{l}{2}}=0,~~l=2,4\cdots, 2M-2,\\
&(\frac{\tau}{2} A+\beta \mathcal{M}_l^T
\mathcal{M}_l)\widetilde{z}_l^k + \beta\mathcal{M}_l^T\bigl(
         q^k-\mathcal{M}_lz_l^k \bigr)-\frac{\tau}{2} d_{\frac{l}{2}}=0,~~l=2M.
\end{aligned}
\end{eqnarray*}
By the definitions of $\mathcal{M}_l$ in \eqref{eq: A-Problem-M}, we can derive
\begin{eqnarray*}
\begin{aligned}
&\left(\tau A+\beta \big
((A+\frac{\tau}{2}B)^T(A+\frac{\tau}{2}B)+(A-\frac{\tau}{2}B)^T(A-\frac{\tau}{2}B)\big)\right)(\widetilde{z}_2^k,
\widetilde{z}_4^k, \cdots, \widetilde{z}_{2M-2}^k)\\
=&\tau
(d_1,d_2,\cdots,d_{M-1})-\beta\left((A+\frac{\tau}{2}B)^T(q_1^k,q_2^k,\cdots,q_{M-1}^k)+(-A+\frac{\tau}{2}B)^T(q_2^k,q_3^k,\cdots,q_{M}^k)\right)\\
&+\beta\left((A+\frac{\tau}{2}B)^T(A+\frac{\tau}{2}B)+(A-\frac{\tau}{2}B)^T(A-\frac{\tau}{2}B)\right)(z_2^k,
z_4^k, \cdots, z_{2M-2}^k),\\
&\left(\frac{\tau}{2} A+\beta (A+\frac{\tau}{2}B)^T(A+\frac{\tau}{2}B)\right)\widetilde{z}_{2M}^k\\
=&\frac{\tau}{2} d_{M}-\beta(A+\frac{\tau}{2}B)^T q_{M}^k+\beta
(A+\frac{\tau}{2}B)^T(A+\frac{\tau}{2}B) z_{2M}^k.
\end{aligned}
\end{eqnarray*}
Let $D_{y,z}=2A^TA+\dfrac{\tau^2}{2}B^TB$,
$D_{q,1}=(A+\dfrac{\tau}{2}B)^T$, $D_{q,2}=(-A+\dfrac{\tau}{2}B)^T$,
and $D_{y,M}=A^TA+\tau A^TB+\dfrac{\tau^2}{4}B^TB$, we obtain
\begin{eqnarray*}
\begin{aligned}
&(\tau A+\beta D_{y,z})(\widetilde{z}_2^k,
\widetilde{z}_4^k, \cdots, \widetilde{z}_{2M-2}^k)\\
=&\tau(d_1,d_2,\cdots,d_{M-1})-\beta[D_{q,1}(q_1^k,q_2^k,\cdots,q_{M-1}^k)+D_{q,2}(q_2^k,q_3^k,\cdots,q_{M}^k)]\\
&+\beta D_{y,z}(z_2^k,
z_4^k, \cdots, z_{2M-2}^k),\\
&(\frac{\tau}{2} A+\beta D_{y,M})\widetilde{z}_{2M}^k =\frac{\tau}{2}
d_{M}-\beta D_{q,1} q_{M}^k+\beta D_{y,M} z_{2M}^k.
\end{aligned}
\end{eqnarray*}
Further more, let $q_{M+1}^k=\mathbf{0}$, and
\begin{eqnarray*}
\kappa=\left\{\begin{array}{rcl} &1, \quad& m=2,\cdots 2M-2,\\[0.2cm]
&\frac{1}{2},  \quad & m=2M,
\end{array}\right.\qquad
D_y=\left\{\begin{array}{rcl} &D_{y,z}, \quad& m=2,\cdots 2M-2,\\[0.2cm]
&D_{y,M},  \quad & m=2M,
\end{array}\right.
\end{eqnarray*}
the first-order optimality conditions for the subproblems for the state variable
can be rewritten as
\begin{eqnarray}\label{ADMMy}
\begin{aligned}
&(\kappa \tau A+\beta D_{y})(\widetilde{z}_2^k,
\widetilde{z}_4^k, \cdots, \widetilde{z}_{2M}^k)\\
=&\kappa\tau(d_1,d_2,\cdots,d_{M})-\beta[D_{q,1}(q_1^k,q_2^k,\cdots,q_{M}^k)+D_{q,2}(q_2^k,q_3^k,\cdots,q_{M+1}^k)]\\
&+\beta D_{y}(z_2^k, z_4^k, \cdots, z_{2M}^k).
\end{aligned}
\end{eqnarray}

It is not difficult to find that the left matrices in
\eqref{ADMMu} and \eqref{ADMMy} are symmetric positive definite, and
the right hands are given matrices. For the column $i$ of the right
hand matrix in \eqref{ADMMu} or \eqref{ADMMy}, we could obtain
$\widetilde{z}_{2i-1}^k$ or $\widetilde{z}_{2i}^k$ directly.
Therefore, the equations \eqref{ADMMu} and \eqref{ADMMy} could be
solved in parallel.

Now, we are at the stage to establish the explicit parallel
implementation of the splitting algorithm 1 as follows.

\textbf{Parallel implementation of splitting algorithm 1}

Input: $\mathcal{A}_m$, $\mathcal{B}_m$, $d_m$,
$(m=1,\cdots,M)$ $\mathcal{F}$, $A$, $D_{u,q}$, $D_{u,z}$, $D_y$,
$D_{q,1}$, $D_{q,2}$, $H$,

\hspace{1.2cm}$\alpha$, $\beta$, $\gamma$, $\nu$, $\tau$,
$\varepsilon$, and initial data
$w^1=(U_{\frac{1}{2}}^1,Y_{1}^1,U_{\frac{3}{2}}^1,\cdots,Y_{M-1}^1,U_{M-\frac{1}{2}}^1,Y_{M}^1,\lambda^1)$.

\begin{itemize}
  \item For $k=1,2,\cdots$

  Let $q^k=\sum\limits_{m=1}^M \mathcal{A}_m Y_m^k+\sum\limits_{m=1}^{M} \mathcal{B}_m U_{m-\frac{1}{2}}^k-\mathcal{F}-\frac{\lambda^k}{\beta}$.

  I$_1$ ~Compute $\widetilde{U}^k$ by the equation
  \eqref{ADMMu} in parallel;

  I$_2$ ~Compute $\widetilde{Y}^k$ by the equation
  \eqref{ADMMy} in parallel;

  I$_3$ ~Compute $\widetilde{\lambda}^k$ by the last equation
  of \eqref{ADMM-GP}.

\hspace{-0.5mm}II$_1$  $w^{k+1}  = w^k- \nu(w^k-\widetilde{w}^k)$;

\hspace{-0.5mm}II$_2$ If $\|w^k-w^{k+1}\|_H^2 \leq\varepsilon$, let $w^*=w^{k+1}$, break.
\item end
\end{itemize}

\begin{rem}
In general, the optimal control problems \eqref{eq:
objection}-\eqref{eq: constraint} always have other limitations on
control variable $u$ or state variable $y$ for the application
purpose. For example, if the state variable has the box constraints
$y_a\leq y \leq y_b$, then the subproblem of (\ref{ADMM-GP}) in $y$
does not have a closed form solution. Fortunately, if we reformulate
it in the following form, our algorithm still works well.

\begin{eqnarray*}
 \begin{array}{rl}
     \min\limits_{Y,U,P} &\left\{ \mathcal{J}_h(Y,U)+\sum\limits_{m=1}^{M} \chi_m(P_m) \right\}   \\[0.2cm]
     s.t.~ & \left[
      \begin{array}{ccc}
        \mathcal{A}& ~\mathcal{B}& \mathbf{0}\\
        \mathbf{I} & ~\mathbf{0} & -\mathbf{I}
        \end{array}\right]
      \left[
      \begin{array}{c}
         Y\\
         U\\
         P
        \end{array}\right]
       =\left[
      \begin{array}{c}
        \mathcal{F}\\
        \mathbf{0}
        \end{array}\right],
        \quad Y,U,P\in \mathbb{R}^{MN_i}~\mbox{or}~\mathbb{R}^{MN}, \quad y_a\leq P \leq
        y_b,
          \end{array}
\end{eqnarray*}
where $\chi_m$ stands for the indicative function on $[y_a, y_b]$,
that is,
\begin{equation*}
\chi_m(y)=
\left
\{ \begin{aligned}
&0,\quad ~~~y\in [y_a, y_b],\\
&+\infty, ~\mbox{else}.
\end{aligned} \right.
\end{equation*}
Now, the above problem could be solved by the splitting algorithm 1
explicitly. Let $(\widetilde{\lambda}^k,\widetilde{\mu}^k)$ be the
lagrangian multipliers at the $k$-th step. By direct calculations,
at the $k$-th step, we obtain

\begin{eqnarray}\label{ADMMupy}\hspace{-7mm}
\begin{aligned}
(\widetilde{U}_{\frac{1}{2}}^k, \widetilde{U}_{\frac{3}{2}}^k,
\cdots, \widetilde{U}_{M-\frac{1}{2}}^k)=~&\beta (\alpha
D_{u,q}+\beta D_{u,z})^{-1}\left(D_{u,z}(U_{\frac{1}{2}}^k,
U_{\frac{3}{2}}^k, \cdots, U_{M-\frac{1}{2}}^k)-D_{u,q}(q_1^k,
q_2^k, \cdots, q_{M}^k)\right),
\\
(\widetilde{P}_1^k,\widetilde{P}_2^k,\cdots,\widetilde{P}_{M}^k)
=~&\mathscr{P}_{[y_a,y_b]}\left((Y_1^k,
Y_2^k, \cdots, Y_{M}^k)-\frac{1}{\beta}(\mu_1^k,\mu_2^k,\cdots,\mu_{M}^k)\right),\\
(\widetilde{Y}_1^k, \widetilde{Y}_2^k, \cdots, \widetilde{Y}_{M}^k)
=~&(\kappa \tau A+\beta D_{y}+\beta I)^{-1}\left( \beta(D_{y}(Y_1^k,
Y_2^k, \cdots, Y_{M}^k)+(P_1^k,P_2^k,\cdots,P_{M}^k)\right.\\
&-D_{q,1}(q_1^k,q_2^k,\cdots,q_{M}^k)-D_{q,2}(q_2^k,q_3^k,\cdots,q_{M+1}^k))\\
&\left.+(\mu_1^k,\mu_2^k,\cdots,\mu_{M}^k)+\kappa\tau(d_1,d_2,\cdots,d_{M})\right),\\
(\widetilde{\lambda}^k,~\widetilde{\mu}^k)=~&\left({\lambda}^k-\beta(\mathcal
{A} \widetilde{Y}^k+\mathcal {B} \widetilde{U}^k-\mathcal {F}),~
{\mu}^k-\beta( \widetilde{Y}^k-\widetilde{P}^k)\right).
\end{aligned}
\end{eqnarray}
where $\mathscr{P}_{[y_a,y_b]}$ denotes the projection operator from
$\mathbb{R}^{MN_i}$ for Dirichlet case(or $\mathbb{R}^{MN}$ for
Neumann case) to $[y_a,y_b]$.
\end{rem}

\section{Numerical experiments}\label{NE}
In this section, we present some numerical examples to verify our
theoretical results in section \ref{ajsa}. Let
$\Omega=[0,1]\times[0, 1]\subset \mathbb {R}^2$, and consider the
parabolic optimal control problem \eqref{eq: objection}-\eqref{eq:
constraint} with Dirichlet boundary condition or Neumann boundary
condition.


In order to verify the convergence results in Theorem \ref{thm: Allconverge},
we need to check the convergence order by the finite element method (FEM)
in Lemma \ref{lem:fem} and the iteration error of full Jacobian
decomposition method with correction in Lemma \ref{lem: contraction}, separately.
Since the convergence rate with respect to the grid size
$h$ and the time step size $\tau$ are all of order two
theoretically, we can choose $M=\frac{T}{h}$, which allows us to use
$\tau=h$. Therefore, we only need to verify the convergence order
with respect to $h$, which is the same as $\frac{1}{\sqrt{N_f}}$,
where $N_f$ is the degree of freedom in spatial direction. We adopt
two strategies to test the convergence of our proposed method: (1)
choose the number of iteration $k$ large enough (e.g. $k=10^4$) and
compute the errors of FEM on nested triangulations with refined mesh
size to see the convergence order of FEM. (2) choose $h$ small
enough (e.g. $h =1/64$) and compute the iteration errors by the
full Jacobian decomposition method with correction to see the convergence rate of
Algorithm~1. The FEM discretization error and the correction full
Jacobian decomposition method iteration error are measured by the
$L^2$-norm and $H$-norm, respectively.

For illustrating the efficiency of the proposed parallel algorithm,
we introduce the concept of the parallel speedup factor(PSF). Let $T_c$ and $T_{g}$ be the
total temporal costs of the algorithm on Central Processing Unit(CPU) with serial execution
and Graphics Processing Unit(GPU) running in parallel, then the PSF can be defined as
\begin{align}\label{PSF}
\mbox{PSF}=\frac{T_c}{T_g}.
\end{align}
At the $I_1$ and $I_2$ steps of the $k-$th iteration in ``\textbf{Parallel implementation of splitting algorithm 1}",
we need to solve a large number of linear systems such as $\mathcal{D}\mathcal{X}_k=\mathcal{R}_k$.
Here $\mathcal{D}$ is the coefficient matrix independent of $k$. $\mathcal{X}_k$ and $\mathcal{R}_k$ are
the solution matrix and right-hand matrix respectively, which are all dependent on $k$. Instead of using the formula
$\mathcal{X}_k=\mathcal{D}\setminus \mathcal{R}_k$ to solve the systems of linear equations on CPU called \textbf{CPU-algorithm},
we will use $\mathcal{X}_k=\mathcal{D}^{-1}\mathcal{R}_k$ to solve it on GPU called \textbf{GPU-algorithm}. One reason is
that we only need to compute the inverse of $\mathcal{D}$ once. Moreover, with the help of large scale parallel operation,
GPU has obvious advantages in dealing with matrix multiplication compared with CPU. In fact, the following examples
all could show GPU-algorithm is much faster than CPU-algorithm.

The initial values in all examples are set to be zeros. And all these simulations
are implemented on a computer with a 2.9GHz CPU named "Intel(R) Xeon(R) Platinum 8268" and a GPU named "TITAN V"
with "ComputeCapabiliy=7.0, MaxThreadsPerBlock=1024" by Matlab.

\noindent {\bf Example 5.1.}~(Dirichlet boundary condition) Assume
$T=2$ in \eqref{eq: objection}-\eqref{eq: constraint}. Let the state
function and source function be\vspace{-3mm}
\begin{eqnarray*}
y_d&=&(\cos (\pi t)-\alpha \pi \cos (\pi t)+2\alpha \pi^2 \sin(\pi
t))\sin(\pi x_1)\sin(\pi x_2),\\
f&=&(2\pi^2 \cos(\pi t)-\pi \sin(\pi t)- \sin(\pi t)),
\end{eqnarray*}
respectively. Set the homogeneous Dirichlet boundary condition as
(\ref{eq: boundary}). This example is taken from \cite{LLX17} and
the exact solution of \eqref{eq: objection}-\eqref{eq: constraint}
is
$$y^*=\cos (\pi t)\sin(\pi x_1)\sin(\pi x_2),\quad u^*=\sin (\pi t)\sin(\pi x_1)\sin(\pi x_2).$$

\begin{figure}[htp]
\centering
\includegraphics[width=6.0cm,height=4.0cm]{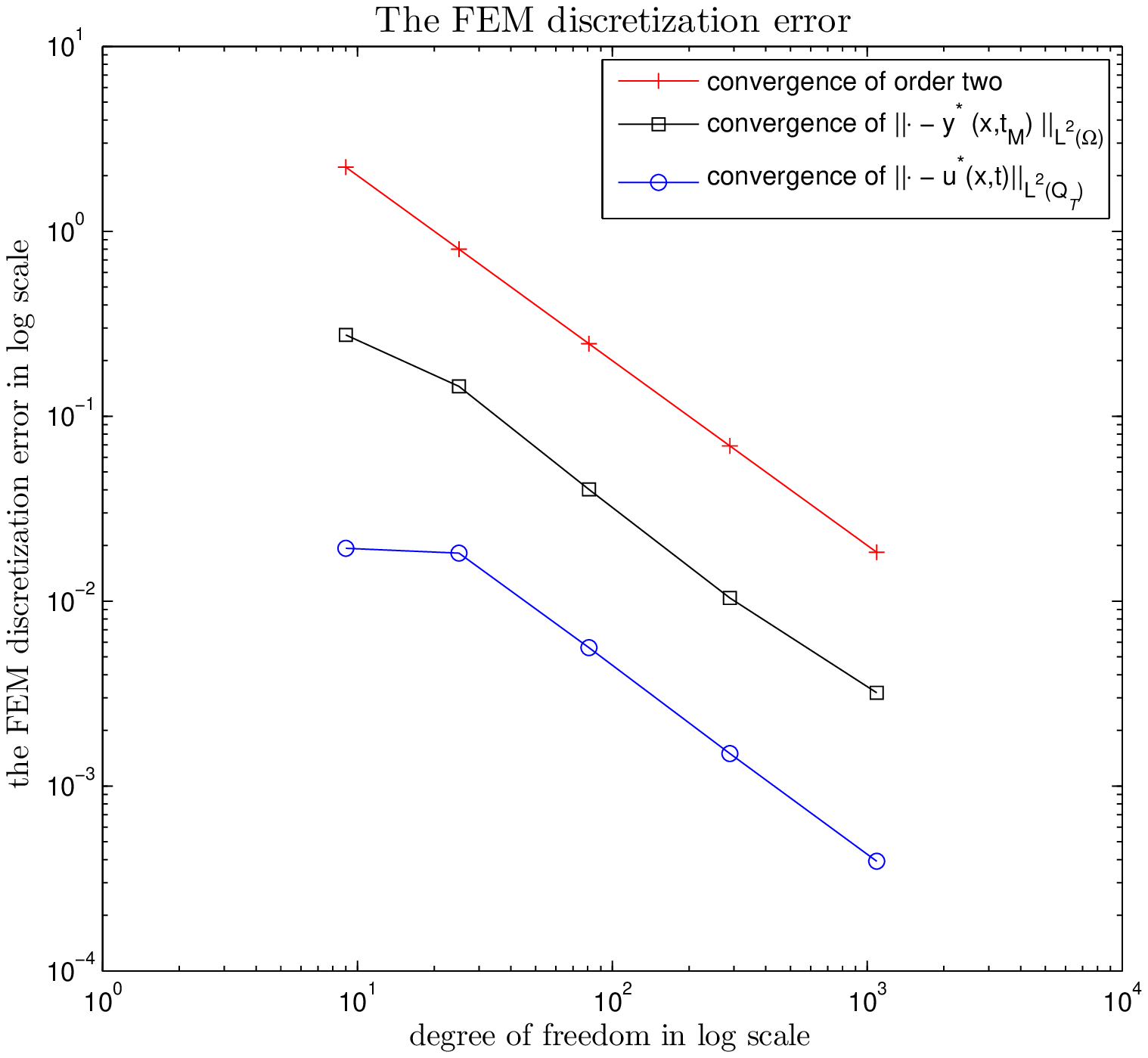}
\includegraphics[width=6.0cm,height=4.0cm]{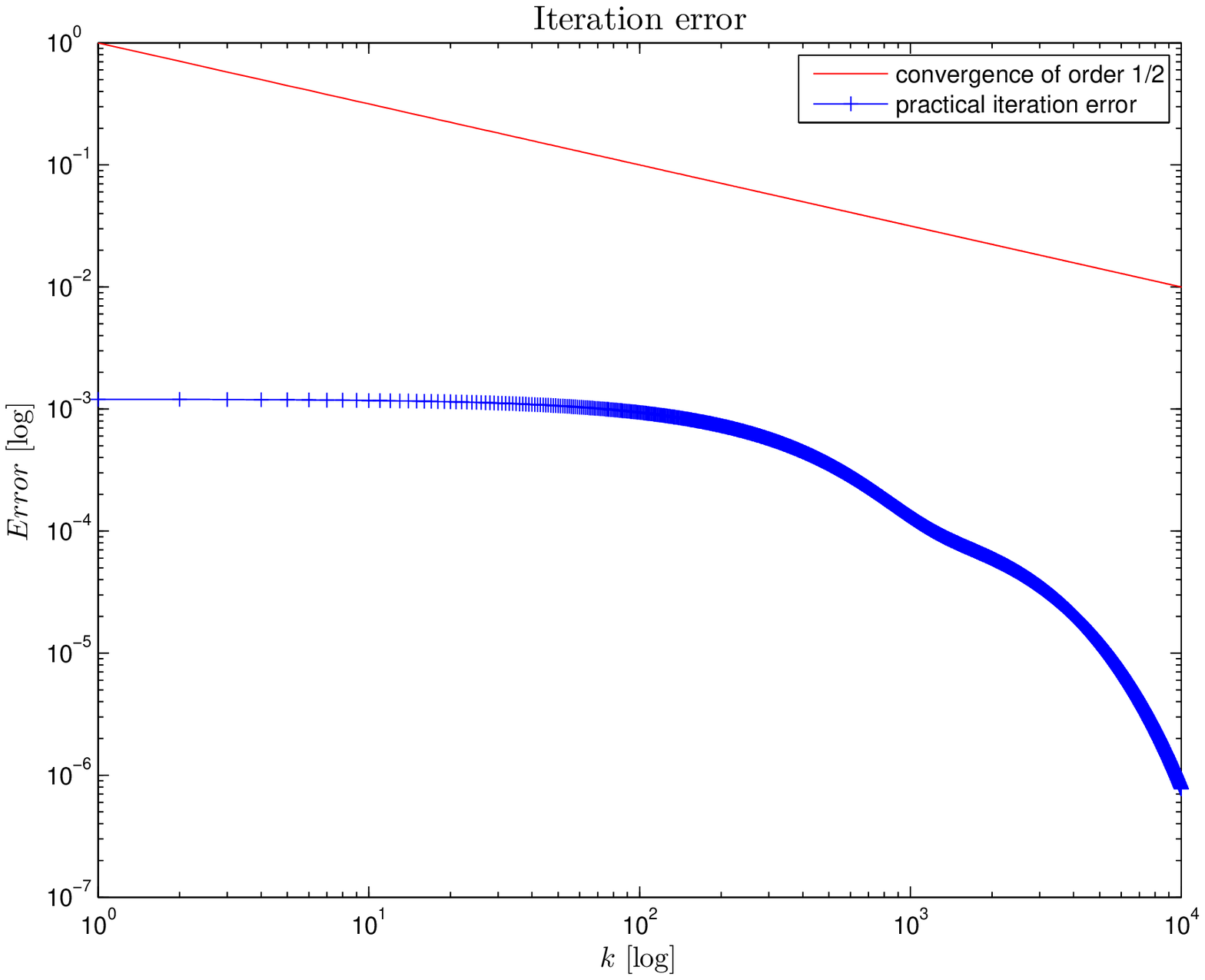}
\caption{\scriptsize Example 5.1, $(\alpha,\beta)=(10^{-2},10)$. The left figure
shows the FEM discretization errors with respect to d.o.f.. The
right figure shows the iterative errors of correction full
Jacobin decomposition method with respect to the $k$-iteration
($\|w^{k}-w^{*}\|_{H}$) in log-scale.} \label{fig:exp11_conv_order}
\end{figure}

Now, we shall test the conclusion in Theorem \ref{thm: Allconverge} by
this example with  $(\alpha,\beta)=(10^{-2},10)$. Figure\,\ref{fig:exp11_conv_order}(left)
shows the convergence results of the FEM as functions of the mesh size $h$ when the
iterative number $k=10^4$ is fixed. We can find that the convergence
rates of the state variable $y$ at $t=T$ with norm
$\|\cdot\|_{L^2(\Omega)}$ and the control variable $u$ with norm
$\|\cdot\|_{L^2(Q_T)}$ are both of order two, which are consistent
with the theoretical results in Theorem \ref{thm: Allconverge}.
Next, we set the mesh size $h=1/32$ as the finest mesh.
Figure\,\ref{fig:exp11_conv_order}(right) shows the errors of
the correction full Jacobian decomposition method ($\|w^{k}-w^{*}\|_{H}$)
with respect to the $k$-iteration in log-scale. It is easy to
see that the full Jacobian decomposition method with correction is faster than
$O(\frac{1}{\sqrt{k}})$ as shown in Lemma~\ref{lem: contraction}.

\begin{figure}[htpb]
\centering
\includegraphics[width=6.0cm,height=4.0cm]{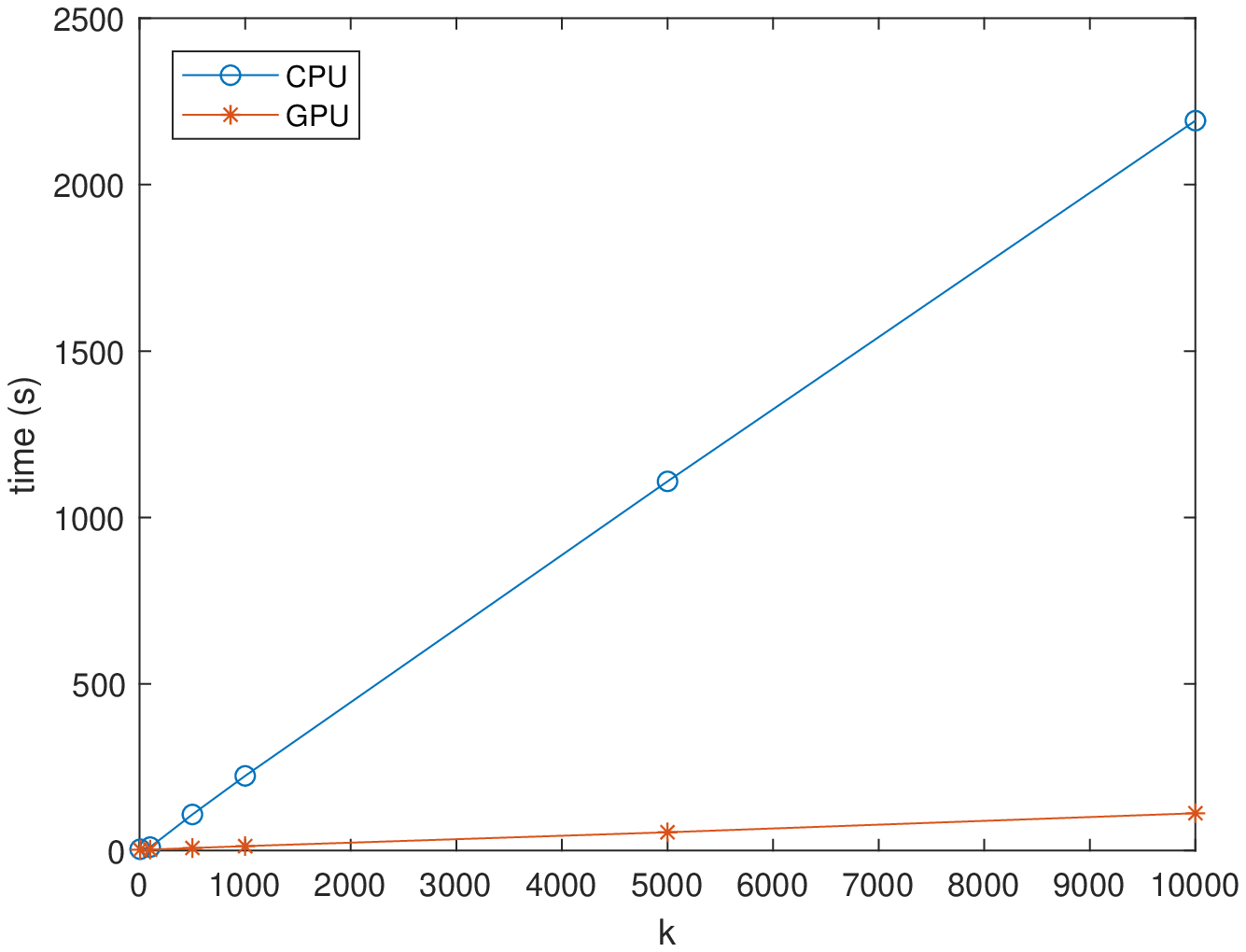}
\includegraphics[width=6.0cm,height=4.0cm]{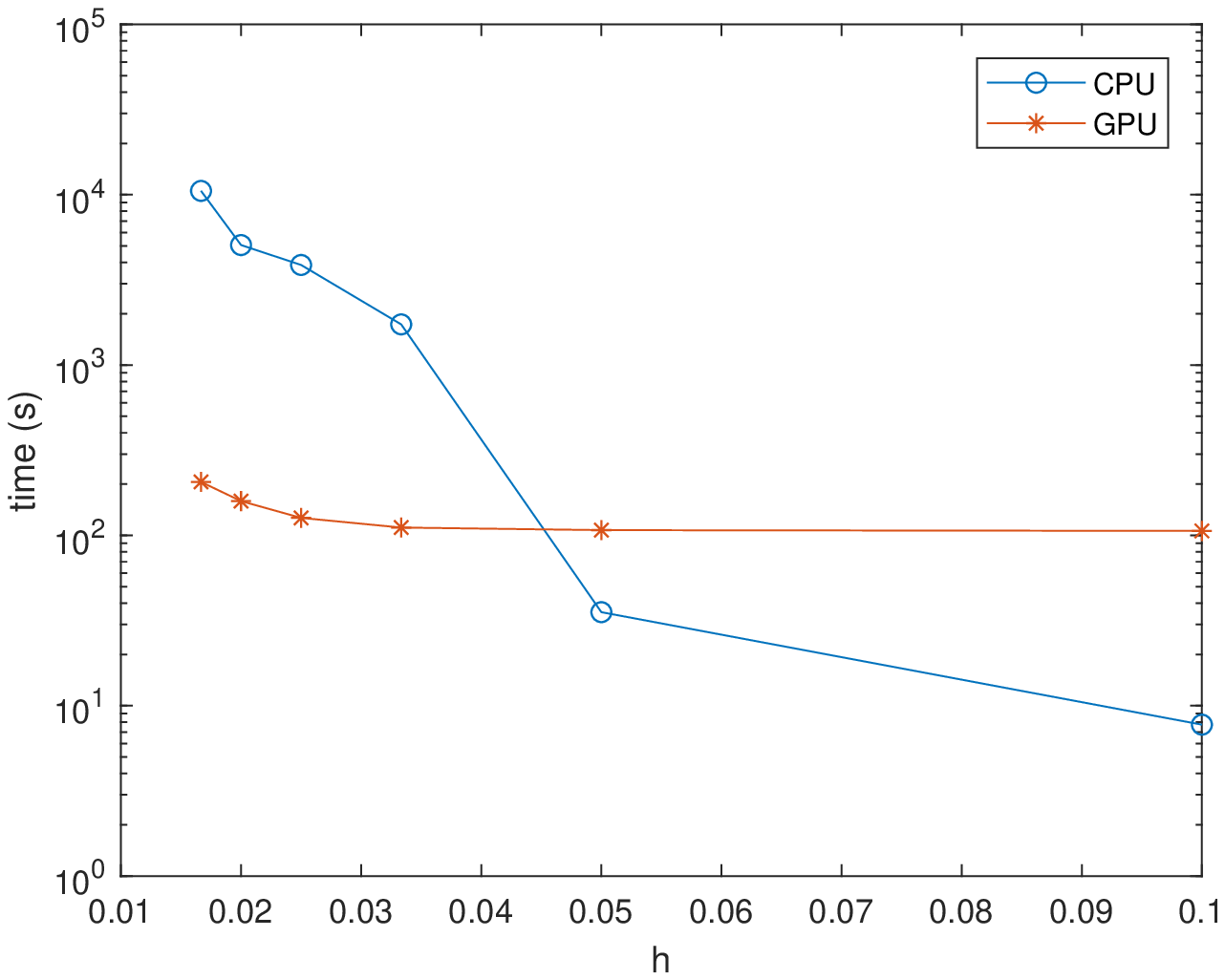}
\caption{\scriptsize Example 5.1, $(\alpha,\beta)=(10^{-2},10)$. The time costs of CPU-algorithm and GPU-algorithm with respect to number of iterations with $h=1/32$ (left) and mesh size with
$k=10^{4}$ (right).} \label{fig:exp11_parallel}
\end{figure}

Further, in order to illustrate the parallel efficiency, we shall consider two cases for Example 5.1:
(I) $h=1/32$ and $k=50,~100,~500,~1000,~5000,~10000$; (II) $k=10^4$ and $h=1/10,~1/20,~1/30,~1/40,~1/50,~1/60$.
 The time costs spent by using CPU-algorithm and GPU-algorithm, and the PSF are recorded
in Table 1. For the sake of intuition, we also present the time costs of Table 1 in Figure~\ref{fig:exp11_parallel}.
From the results in Table 1 and Figure~\ref{fig:exp11_parallel}, we can conclude that the growth rate of time cost
of GPU-algorithm with the number of iterations is much smaller than that of CPU-algorithm in the case of $h=1/32$.
For the fixed iteration number $k=10000$, the time cost of GPU-algorithm changes little with the change of mesh size,
while the time cost of CPU-algorithm changes greatly. Specially, GPU-algorithm reduces the computation time sharply
and PSF reaches $51.12$ when the mesh size $h=1/60$.\vspace{3mm}

{\small
\begin{center}
\begin{tabular}{|c|c|c|c|c|c|c|}\hline
$h=1/32,k=$  &   50   & 100    &   500   & 1000    & 5000     & 10000 \\\hline
CPU(s)       &  3.60  & 9.95   &  108.04 & 223.92  & 1108.62  & 2192.23\\ \hline
GPU(s)       &  2.33  & 2.52   &  7.19   & 12.91   & 54.48    & 111.69\\ \hline
PSF          &  1.55  & 3.95   &  15.03  & 17.34   & 20.35    & 19.63\\ \hline\hline
$k=10^4, h=$ & 1/10   &1/20    & 1/30    &1/40     &1/50      &1/60\\\hline
CPU(s)       & 7.74   &35.41   & 1732.33 &3867.60  &5061.27   &10525.36\\ \hline
GPU(s)       & 106.18 &107.33  & 111.12  &126.65   &158.90    &205.91\\ \hline
PSF          & 0.07   &0.33    & 15.59   &30.54    &31.85     &51.12\\ \hline
\end{tabular}
\vspace{2mm}

{\scriptsize Table~1. Example 5.1, Comparison of CPU-algorithm and GPU-algorithm in computing time (s) and PSF.}
\end{center}
}\vspace{2mm}

\begin{figure}[htpb]
\centering
\includegraphics[width=6.0cm,height=4.0cm]{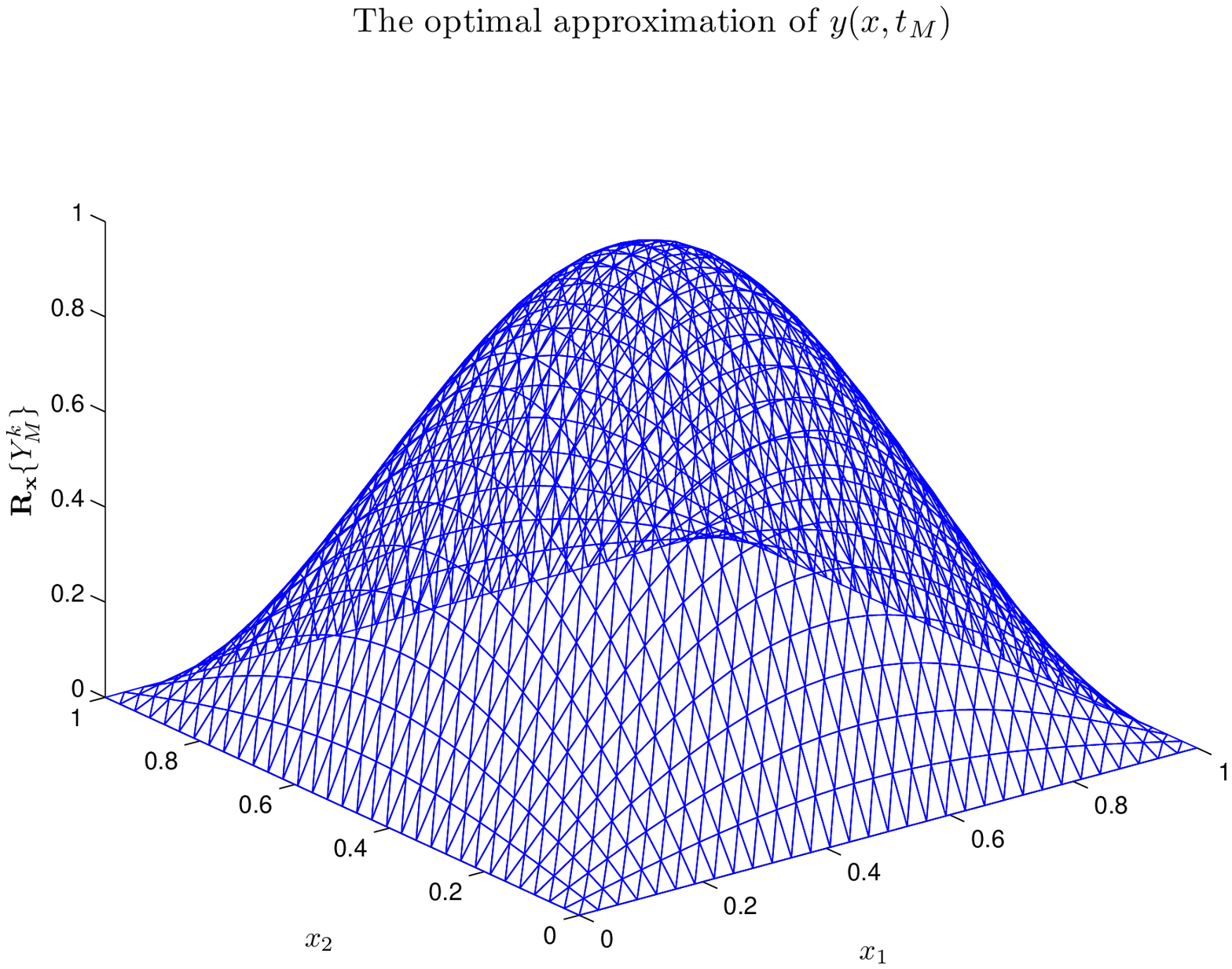}
\includegraphics[width=6.0cm,height=4.0cm]{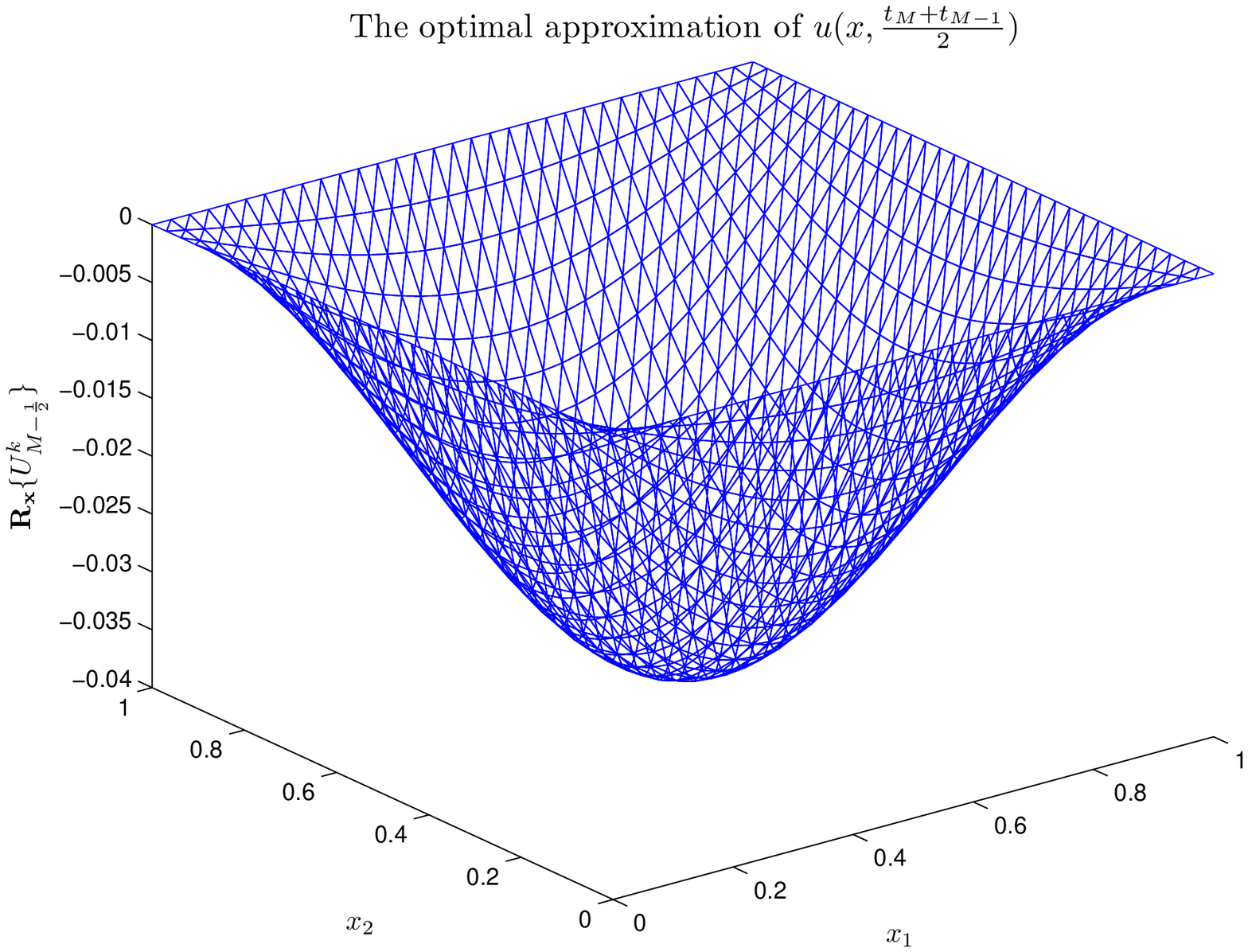}
\caption{\scriptsize Example 5.1, $(\alpha,\beta)=(10^{-2},10)$. The numerical
solutions $\mathbf{R_x}\{Y_{M}^k\}$ and
$\mathbf{R_x}\{U_{M-\frac{1}{2}}^k\}$ with $h=1/32$ and $k=10^{4}$.}
\label{fig:exp11_numerical_solutions}
\end{figure}

For giving the interested readers a visual understanding, we present the numerical
solutions $\mathbf{R_x}\{Y_{M}^k\}$ and $\mathbf{R_x}\{U_{M-\frac{1}{2}}^k\}$ with $h=1/32$ and $k=10^{4}$
in Figure~\ref{fig:exp11_numerical_solutions}. The results illustrate that the numerical solutions are
coincided with the theoretical solutions $y^*$ and $u^*$, which verified the validity of the proposed
method intuitively.

\begin{figure}[htpb]
\centering
\includegraphics[width=6.0cm,height=4.0cm]{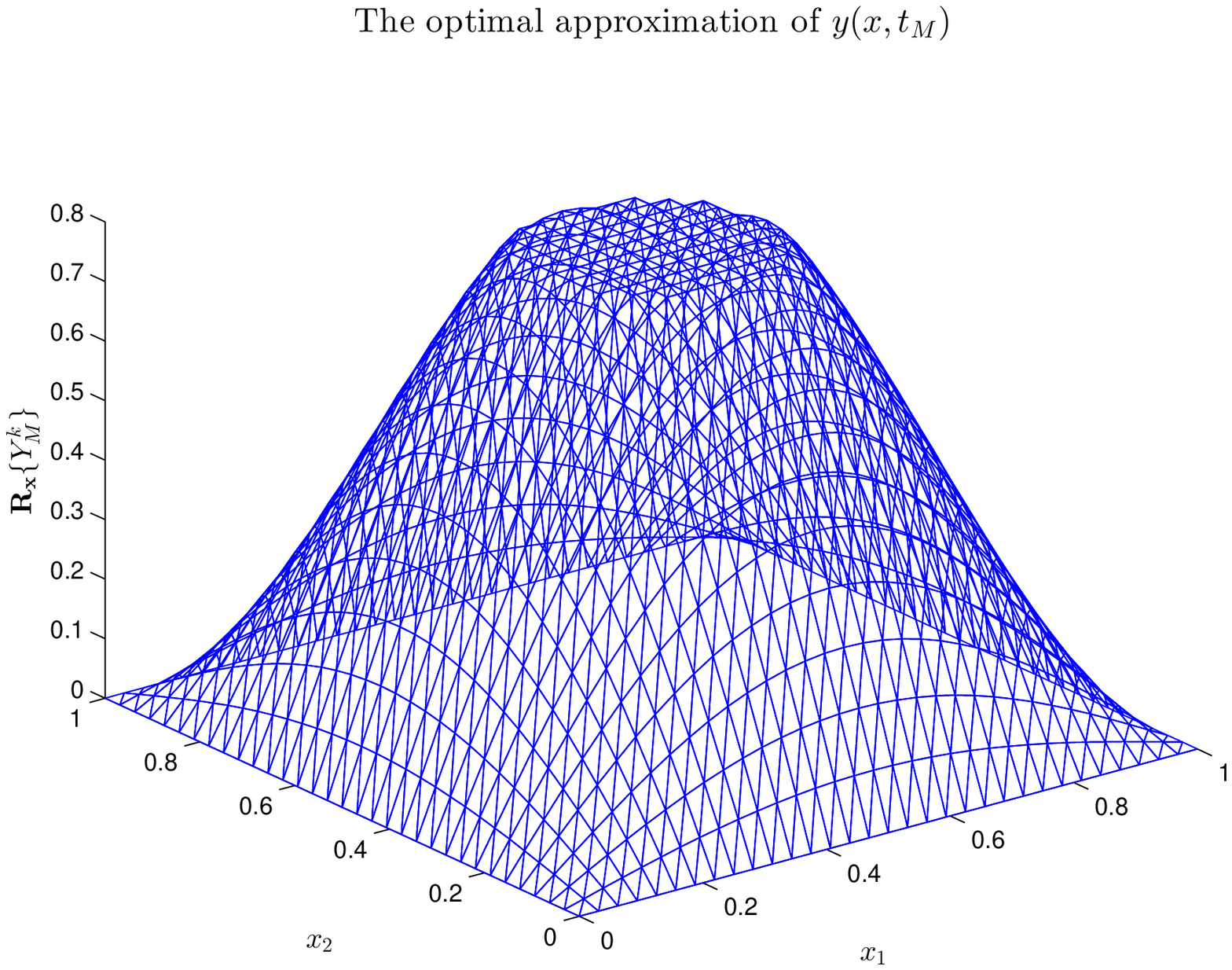}
\includegraphics[width=6.0cm,height=4.0cm]{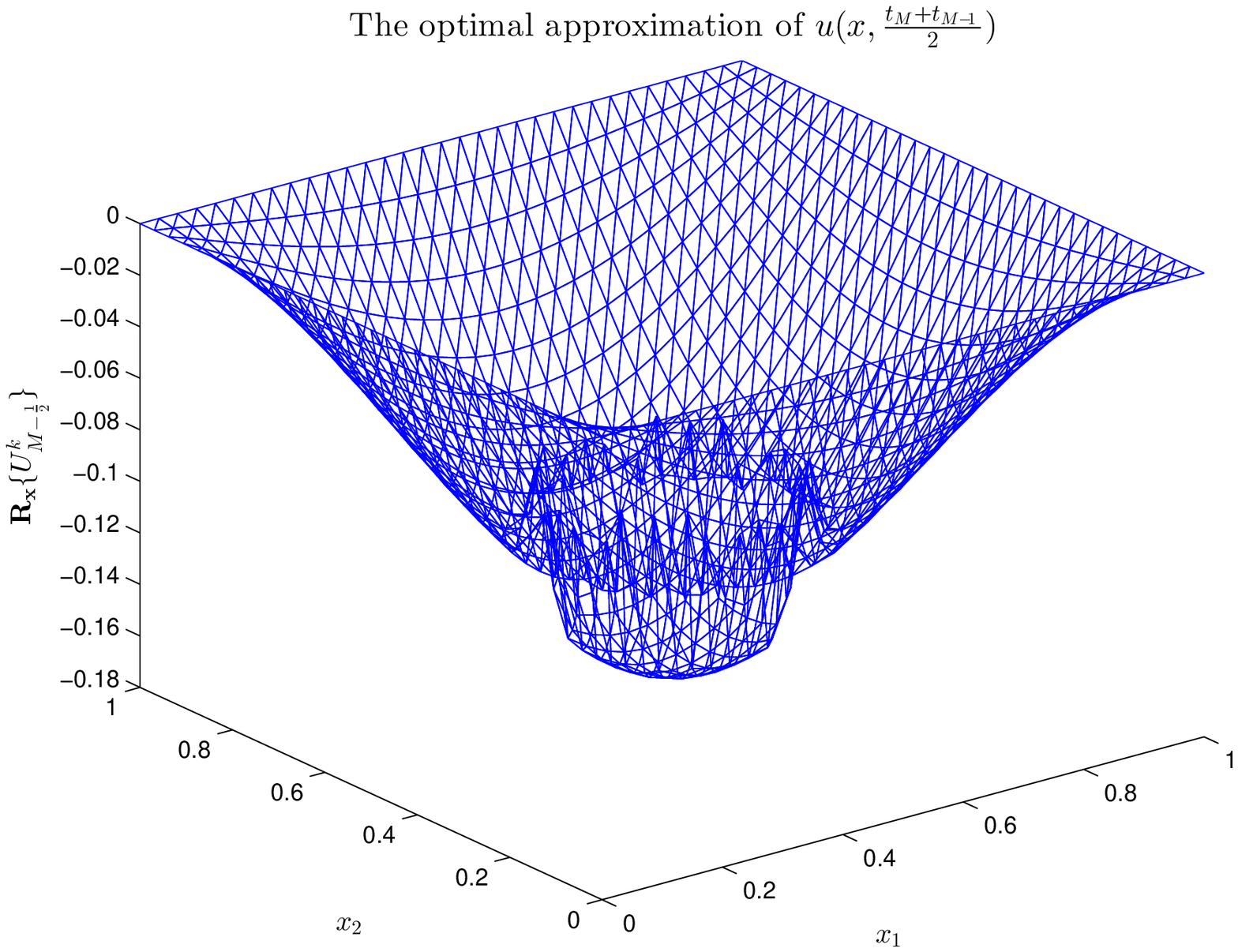}
\caption{\scriptsize Example 5.1, $(\alpha,\beta)=(10^{-2},10)$. The numerical
solutions $\mathbf{R_x}\{Y_{M}^k\}$ and
$\mathbf{R_x}\{U_{M-\frac{1}{2}}^k\}$ with $h=1/32$ and $k=10^{4}$ under the box constraint case.}
\label{fig:exp11_numerical_solutions_box}
\end{figure}

At the end of this example, we carry out the numerical
experiment with the box constraint case $y\in [0,0.8]$. From the
results exhibited in Figure
~\ref{fig:exp11_numerical_solutions_box}, we can see that the
numerical solution is reasonable.

\noindent {\bf Example 5.2.}~(Neumann boundary condition) This example
is taken from \cite{AF12}. Assume $T=1$ in \eqref{eq:
objection}-\eqref{eq: constraint}. Let the expected state function
in the objective function be
$$y_d=c_7\omega_a(x,t)+c_8\omega_b(x,t)+c_9\omega_a(x,0)+c_{10}\omega_b(x,0)+c_{11}\omega_a(x,T)+c_{12}\omega_b(x,T)$$
with
\begin{eqnarray*}
\omega_a(x,t)=e^{\frac{1}{3}\pi^2t}\cos(\pi x_1)\cos(\pi x_2),\quad
\omega_b(x,t)=e^{-\frac{1}{3}\pi^2t}\cos(\pi x_1)\cos(\pi x_2).
\end{eqnarray*}
Set the source function $f=0$ and the homogeneous Neumann boundary
condition as (\ref{eq: boundary}). It follows from the first order
optimality condition \eqref{eq: optcondition} that the exact solution
of \eqref{eq: objection}-\eqref{eq: constraint} is
\begin{eqnarray*}
\begin{aligned}
y^*=&~c_1\omega_a(x,t)+c_2\omega_b(x,t)+c_3\omega_a(x,0)+c_{4}\omega_b(x,0)+c_{5}\omega_a(x,T)+c_{6}\omega_b(x,T),\\
u^*=&~\frac{1}{3}\pi^2(c_1\omega_a(x,t)-c_2\omega_b(x,t))+y^*.
\end{aligned}
\end{eqnarray*}
Here the coefficients $\{c_i\}_{i=1}^{12}$ are specified in Table 2.\vspace{3mm}

{\small
\begin{center}
\begin{tabular}{cccc}\hline
     $c_1$              &     $c_2$          & $c_3$       &  $c_4=c_5=c_6$ \\ \hline
     $\frac{-5\left(5e^{-\frac{1}{3}\pi^2}-6\right)}{-6+7e^{\frac{1}{3}\pi^2}}$
     & $5$
     &$\frac{7+141e^{\frac{1}{3}\pi^2}+7\left(e^{\frac{1}{3}\pi^2}\right)^2-6-106e^{-\frac{1}{3}\pi^2}}{4\left(6-7e^{\frac{1}{3}\pi^2}\right)}$
     & $\frac{1}{4}$\\ \hline
     $c_7$              &     $c_8$          & $c_9$       &  $c_{10}=c_{11}=c_{12}$ \\ \hline
     $\frac{5(9+35 \alpha \pi^4)\left(5e^{-\frac{1}{3}\pi^2}-6\right)}{9\left(6-7e^{\frac{1}{3}\pi^2}\right)}$
     & $5+\frac{175}{9}\alpha \pi^4$
     &$4\cdot c_3 \cdot c_{10}$
     & $\frac{1}{4}+\alpha \pi^4$\\ \hline
\end{tabular}
\vspace{2mm}

{\scriptsize Table~2. The coefficients $\{c_i\}_{i=1}^{12}$ of $y^*$,
$u^*$, and $y_d$.}
\end{center}
}\vspace{2mm}

\begin{figure}[htpb]
\centering
\includegraphics[width=6.0cm,height=4.0cm]{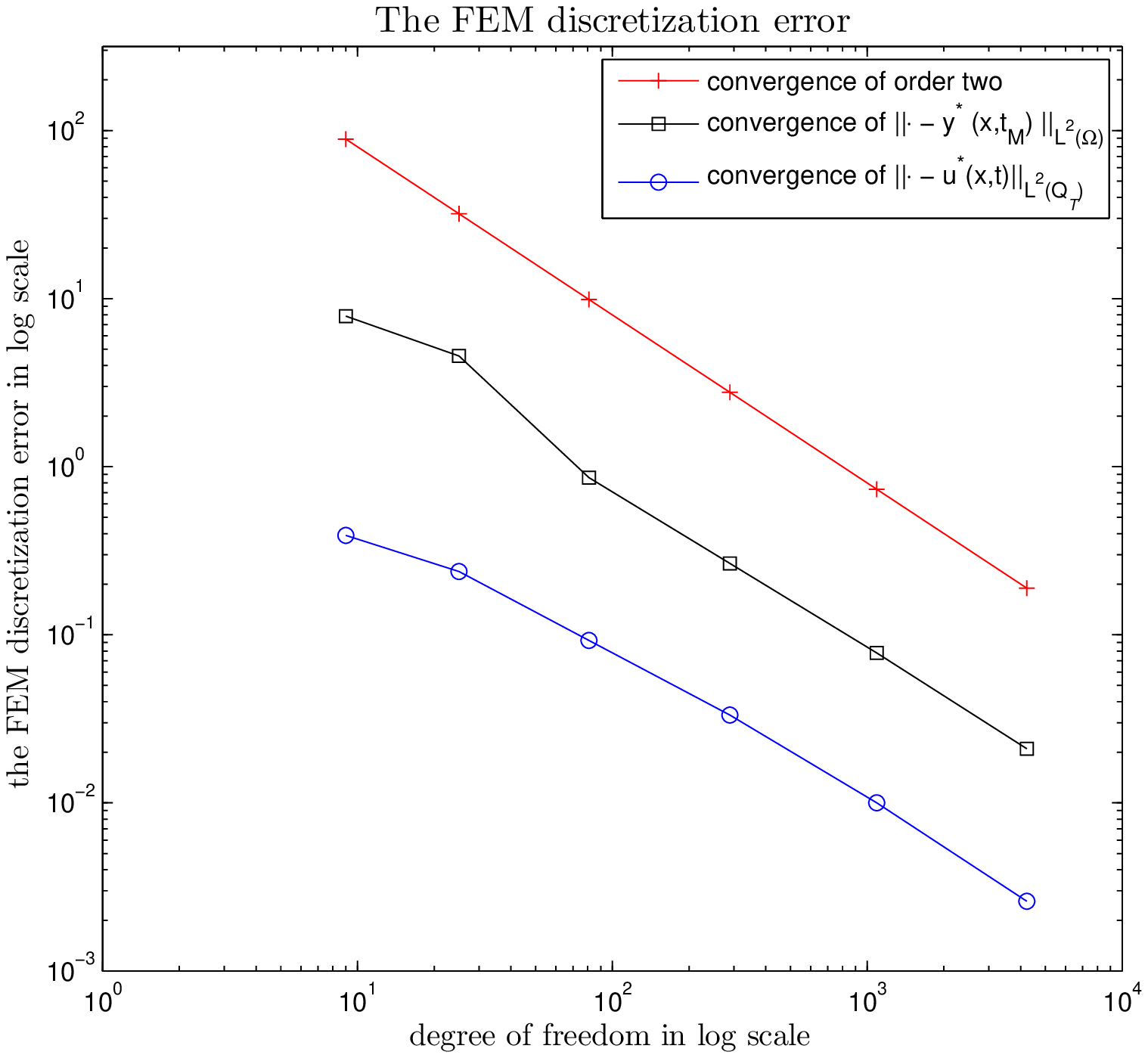}
\includegraphics[width=6.0cm,height=4.0cm]{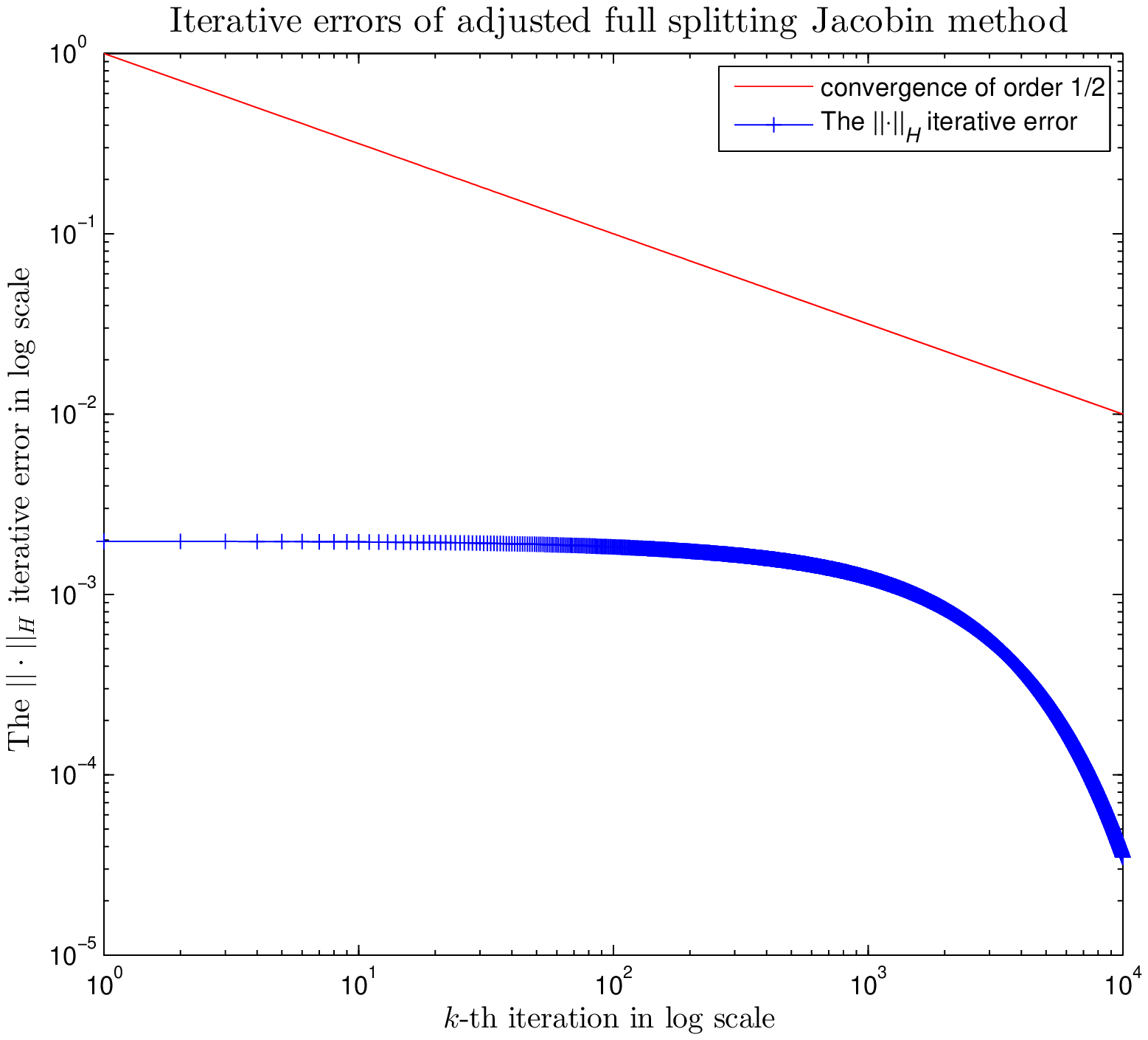}
\caption{\scriptsize Example 5.2, $(\alpha,\beta)=(10^{-3},10^2)$. The left
figure shows the FEM discretization errors with respect to d.o.f..
The right figure shows the iterative errors of correction full
 Jacobin decomposition method with respect to the $k$-iteration
($\|w^{k}-w^{*}\|_{H}$) in log-scale.} \label{fig:exp22_conv_order}
\end{figure}

Now, we verify the convergence result in Theorem \ref{thm: Allconverge} by
Example 5.2 with the parameters $(\alpha,\beta)=(10^{-3},10^2)$.
Figure\,\ref{fig:exp22_conv_order}(left) shows the convergence
results of the FEM as functions of the mesh size $h$ when the
iterative number $k=2\times 10^{4}$ is fixed. We can find that the
convergence rates of the state variable $y$ at $t=T$ with norm
$\|\cdot\|_{L^2(\Omega)}$ and the control variable $u$ with norm
$\|\cdot\|_{L^2(Q_T)}$ are both of order two, which are consistent
with the theoretical results in Theorem \ref{thm: Allconverge}.
Next, we set the mesh size $h=1/64$ as the finest mesh.
Figure\,\ref{fig:exp22_conv_order}(right) shows the errors of
correction full Jacobin decomposition method ($\|w^{k}-w^{*}\|_{H}$)
with respect to the $k$-iteration in log-scale. It is easy to
see that the full Jacobin decomposition method with correction is faster than
$O(\frac{1}{\sqrt{k}})$ as shown in Lemma~\ref{lem: contraction}.

\begin{figure}[htpb]
\centering
\includegraphics[width=6.0cm,height=3.9cm]{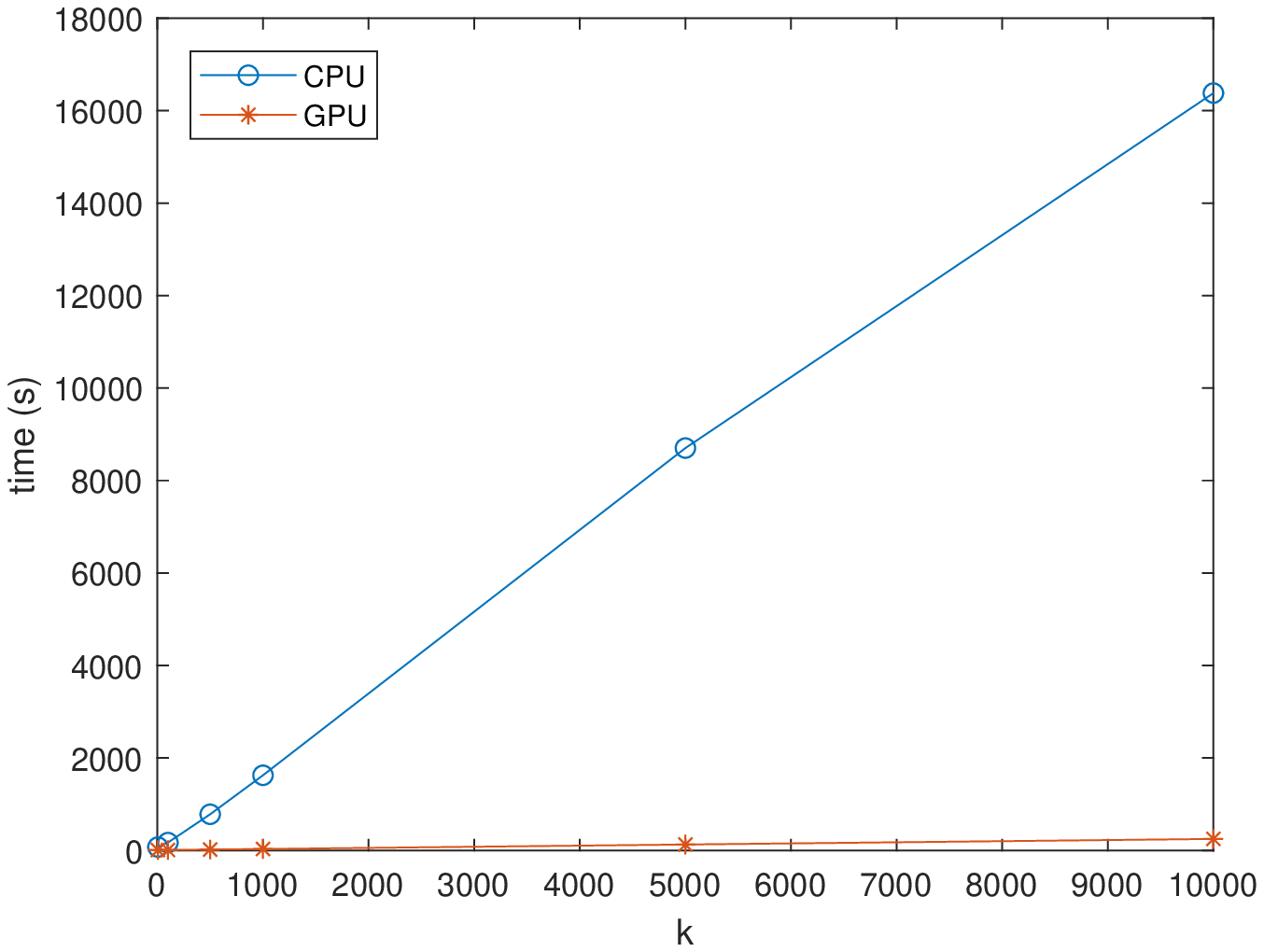}
\includegraphics[width=6.0cm,height=3.9cm]{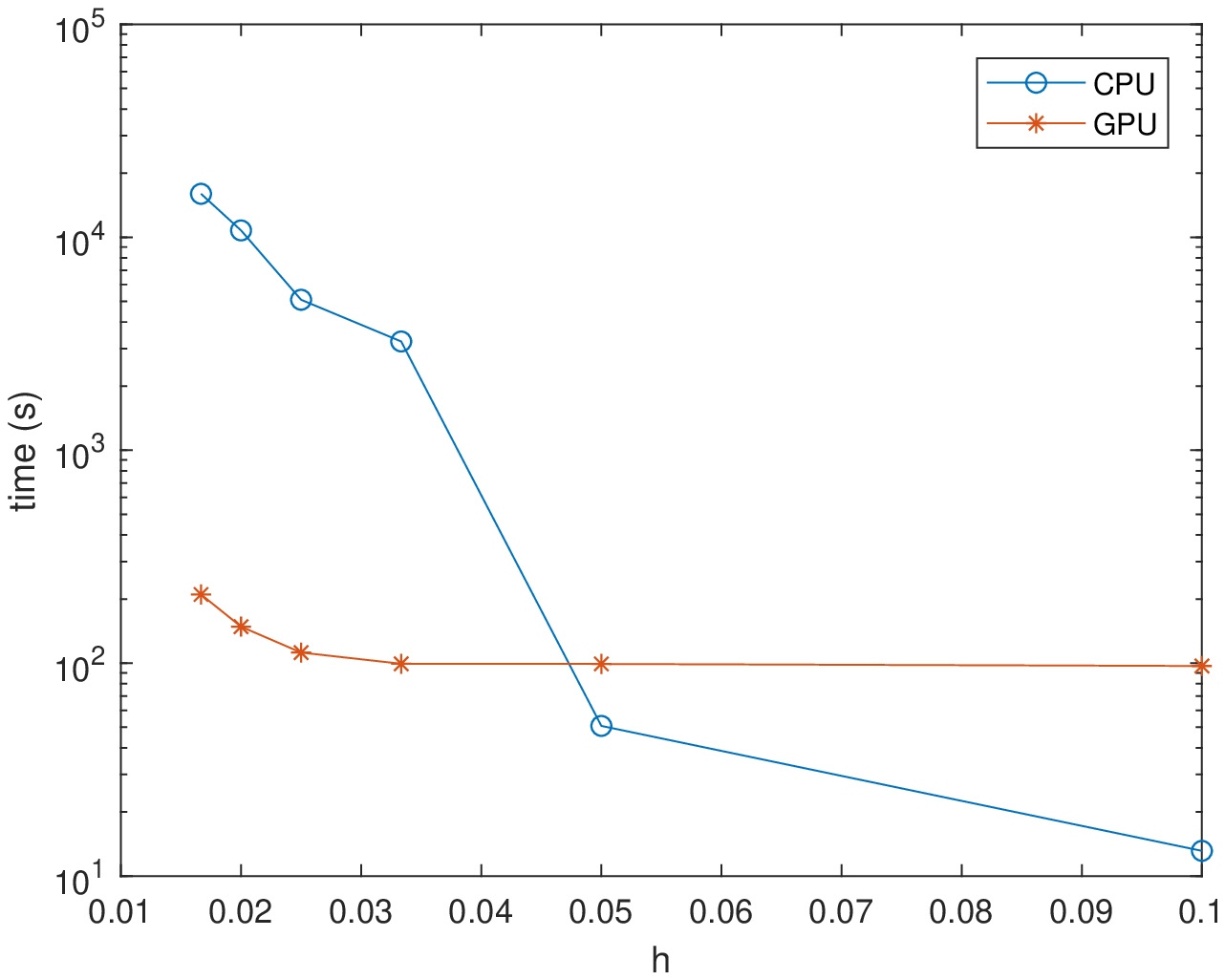}
\caption{\scriptsize Example 5.2, $(\alpha,\beta)=(10^{-3},10^2)$. The time costs of CPU-algorithm and GPU-algorithm with respect to number of iterations with $h=1/64$ (left) and mesh size with $k=10^{4}$ (right).} \label{fig:exp22_parallel}
\end{figure}

Further, the computing times and PSF of Example 5.2 by using CPU-algorithm and GPU-algorithm are listed
in Table 3 and Figure~\ref{fig:exp22_parallel} with (I) $h=1/64$, $k=50,~100,~500,~1000,~5000,~10000$;
(II) $k=10^4$ and $h=1/10,~1/20,~1/30,~1/40,~1/50,~1/60$. Compared with Example 5.1, the problem scale of Example 5.2 is larger,
so the GPU-algorithm is more advantageous and the PSF could reach $76.17$.\vspace{3mm}

{\small
\begin{center}
\begin{tabular}{|c|c|c|c|c|c|c|}\hline
$h=1/64,k=$  &   50   & 100    &   500   & 1000    & 5000    & 10000 \\\hline
CPU(s)       &  73.04 & 169.18 &  781.94 & 1625.20 & 8702.32 & 16380.33\\ \hline
GPU(s)       &  8.45  & 9.97   &  18.91  & 30.74   & 127.73  & 248.44\\ \hline
PSF          &  8.64  & 16.97  &  41.35  & 52.86   & 68.13   & 65.93\\ \hline\hline
$k=10^4, h=$ & 1/10   &1/20    & 1/30    &1/40     &1/50     &1/60\\\hline
CPU(s)       & 13.13  &50.71   & 3240.58 &5087.63  &10769.31 &16001.93\\ \hline
GPU(s)       & 96.95  &99.12   & 99.31   &112.09   &148.35   &210.08\\ \hline
PSF          & 0.13   &0.51    & 32.63   &45.39    &72.59    &76.17\\ \hline
\end{tabular}
\vspace{2mm}

{\scriptsize Table~3. Example 5.2, Comparsion of CPU-algorithm and GPU-algorithm in computing time (s) and PSF.}
\end{center}
}\vspace{2mm}

\begin{figure}[htpb]
\centering
\includegraphics[width=6.0cm,height=3.9cm]{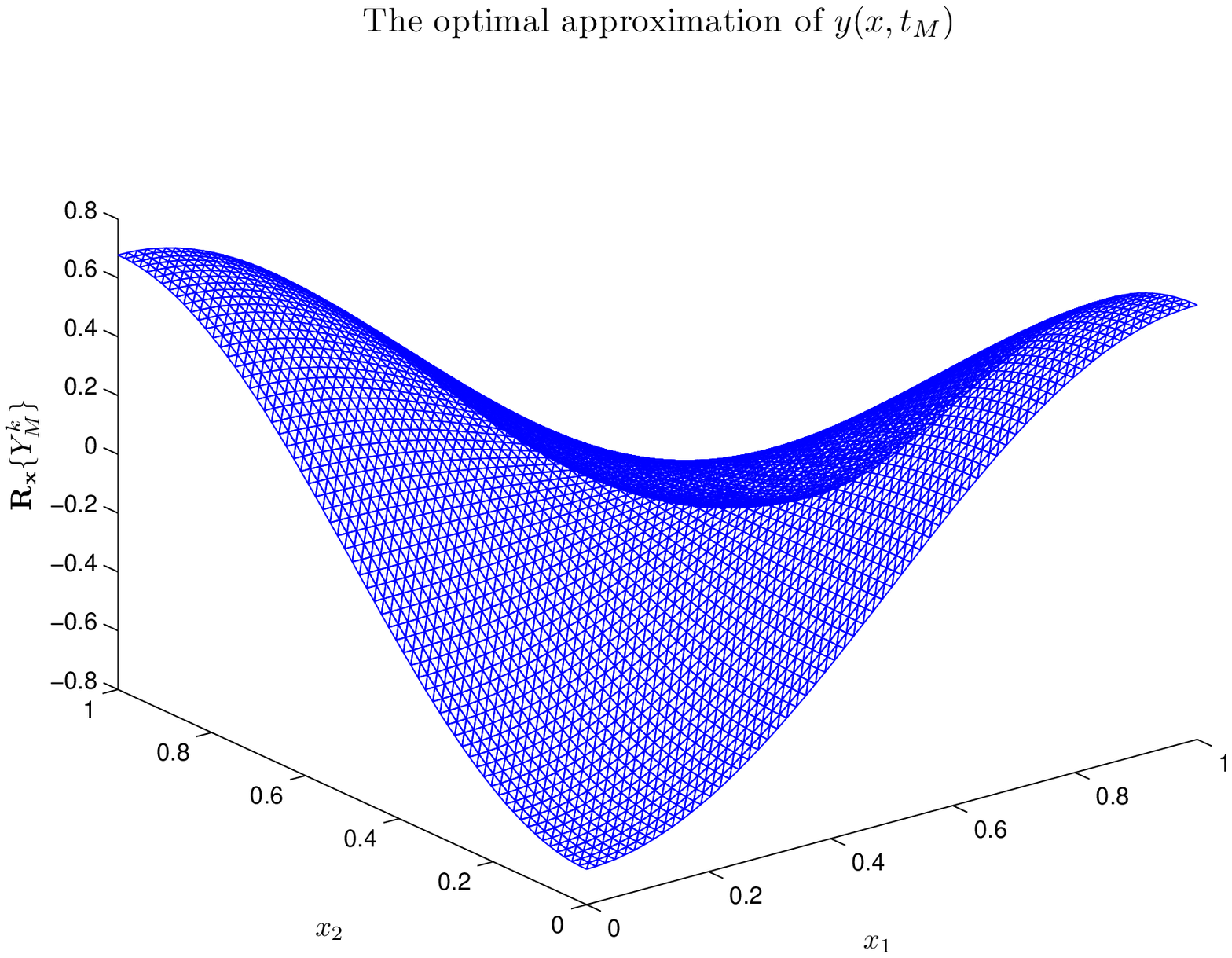}
\includegraphics[width=6.0cm,height=3.9cm]{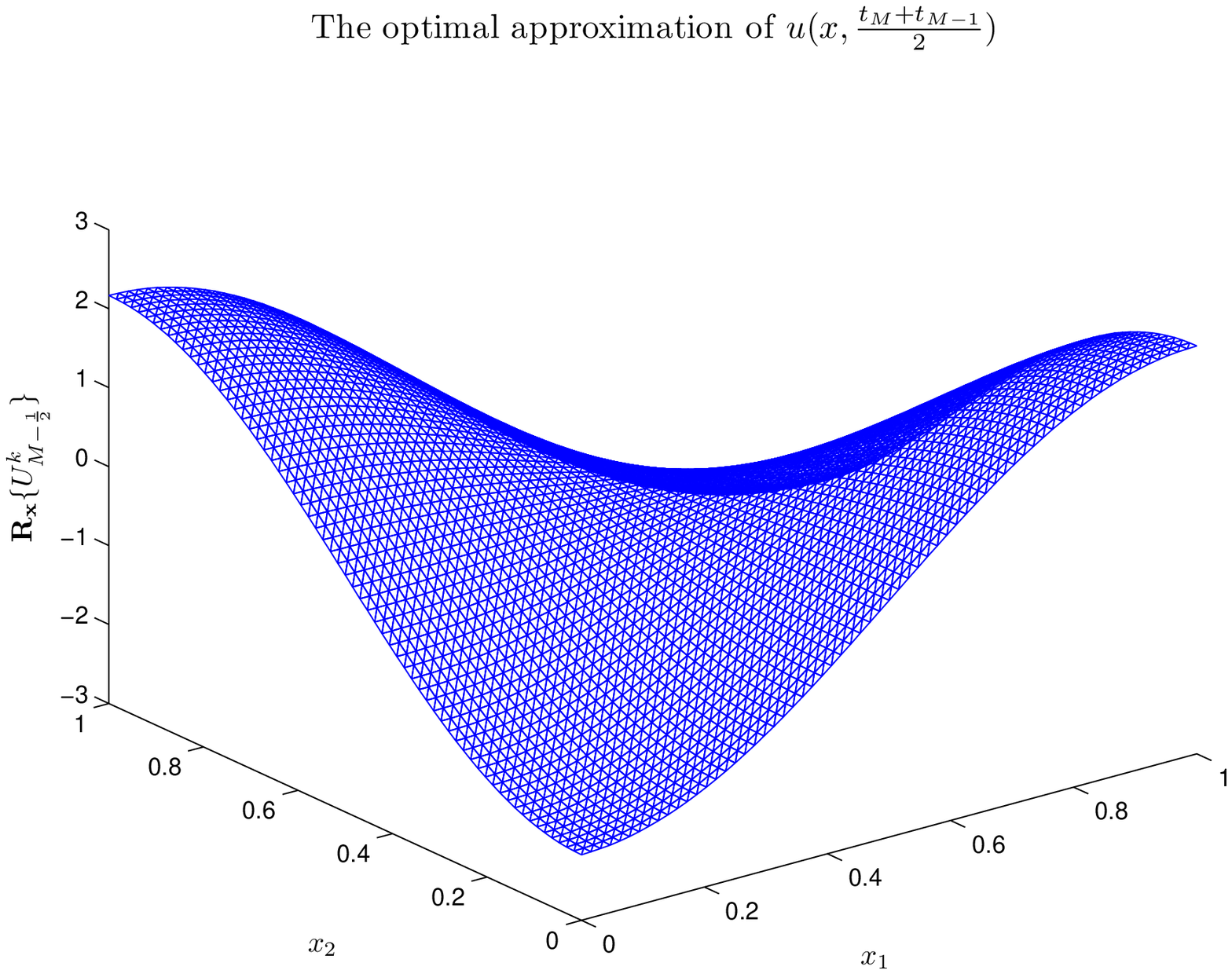}
\caption{\scriptsize Example 5.2, $(\alpha,\beta)=(10^{-3},10^2)$. The numerical
solutions $\mathbf{R_x}\{Y_{M}^k\}$ and
$\mathbf{R_x}\{U_{M-\frac{1}{2}}^k\}$ with $h=1/64$ and
$k=2\times10^{4}$.} \label{fig:exp22_numerical_solutions}
\end{figure}

For giving the interested reader a visual understanding, we present the numerical
solutions $\mathbf{R_x}\{Y_{M}^k\}$ and $\mathbf{R_x}\{U_{M-\frac{1}{2}}^k\}$ with $h=1/64$ and $k=2\times10^{4}$
in Figure~\ref{fig:exp22_numerical_solutions}. The results could illustrate
that the numerical solutions approximate the theoretical solutions
$y^*$ and $u^*$ very well, which verified the validity of the proposed method intuitively.

\begin{figure}[htpb]
\centering
\includegraphics[width=6.0cm,height=3.9cm]{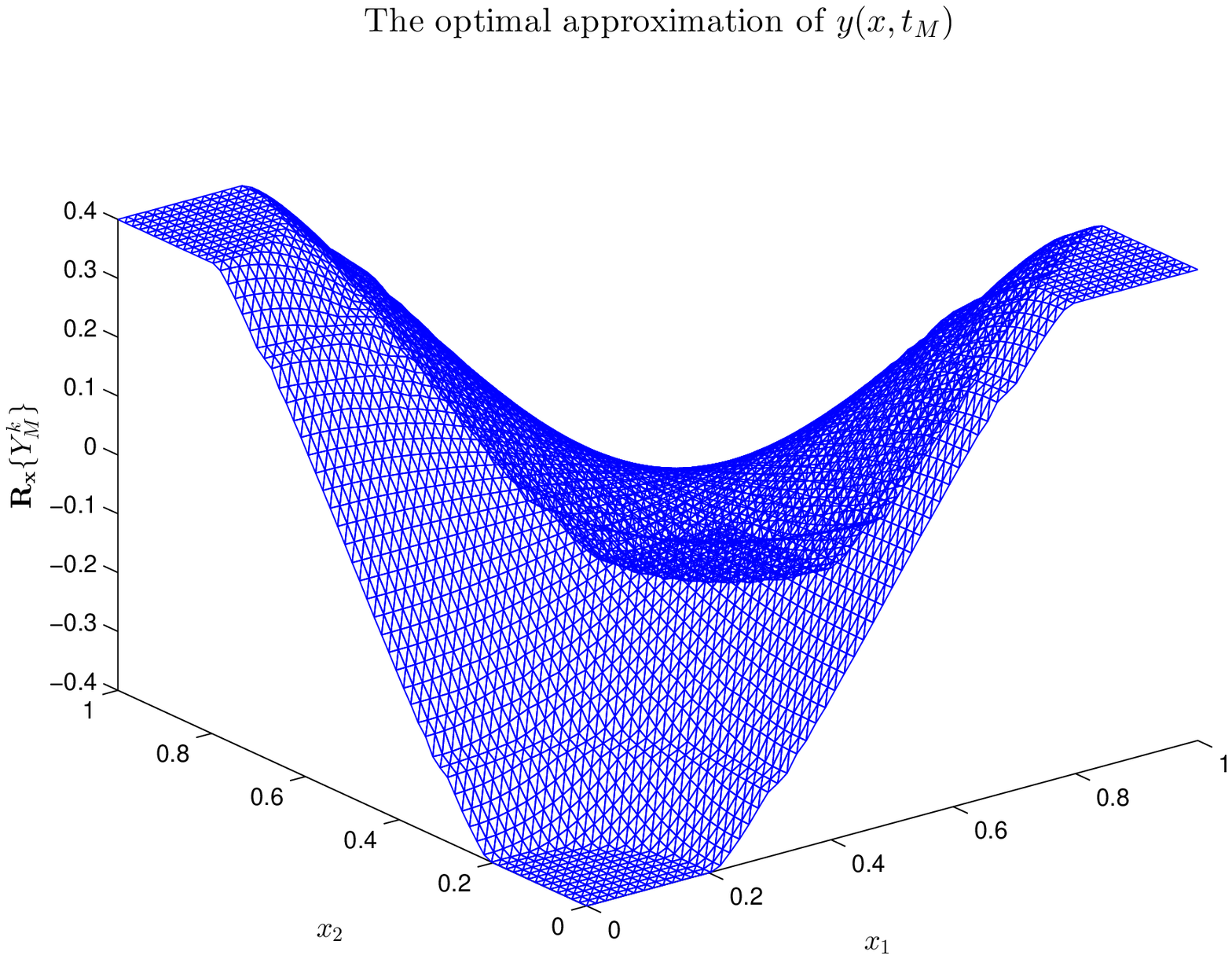}
\includegraphics[width=6.0cm,height=3.9cm]{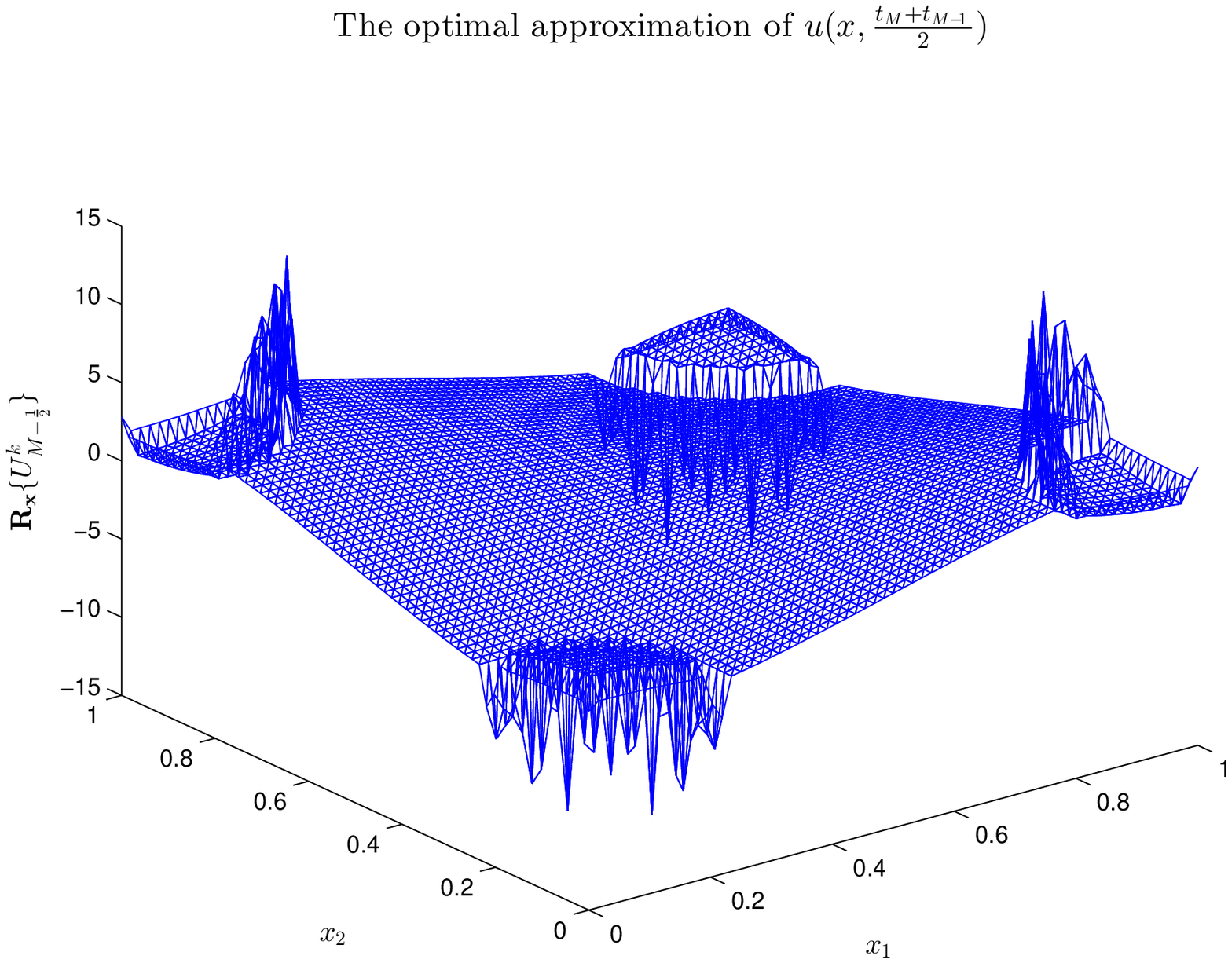}
\caption{\scriptsize Example 5.2, $(\alpha,\beta)=(10^{-3},10^2)$. The numerical
solutions $\mathbf{R_x}\{Y_{M}^k\}$ and
$\mathbf{R_x}\{U_{M-\frac{1}{2}}^k\}$ with $h=1/64$ and
$k=2\times10^{4}$ under the box constraint case.} \label{fig:exp22_numerical_solutions_box}
\end{figure}

Similar to Example 5.1,  we present the numerical
solutions of Example 5.2 with the box constraint case $y\in
[-0.4,0.4]$ in Figure~\ref{fig:exp22_numerical_solutions_box} at the
end of this subsection.

\section{Conclusions}\label{Con}
In this paper, we propose an efficient parallel splitting method for
the parabolic optimal control problems. The model problem is discretized by
the Crank-Nicolson scheme and the numerical integration formula
in temporal direction, and the linear finite element method in
spatial direction.  Based on the separable
structure of the resulting large-scale optimization system, a
full Jacobian decomposition method with correction is proposed, which
improve the computational efficiency significantly. The global
convergence estimation is established based on the FEM discretization
error and the iteration error. Finally, numerical simulations are
presented to verify the efficiency of the proposed algorithm.

\section*{Acknowledgments}
The work of H. Song was supported by the NSF of China under the grant
No. 11701210, the NSF of Jilin Province under the grants No. 20190103029JH, 20200201269JC,
the education department project of Jilin Province under the grant No. JJKH20211031KJ,
and the fundamental research funds for the Central Universities.
The work of J.C. Zhang was supported by the Natural Science Foundation of Jiangsu Province (Grant BK20210540)
, the Natural Science Foundation of the Jiangsu Higher Education Institutions of China
(No. 21KJB110015, 21KJB110001) and the Startup Foundation for Introducing Talent of NJTech (No. 39804131).
The work of Y.L. Hao was supported
by the NSF of China under the grant No. 11901606. The authors
also wish to thank the High Performance Computing Center of Jilin
University, Computing Center of Jilin Province, and Key Laboratory
of Symbolic Computation and Knowledge Engineering of Ministry of
Education for essential computing support.

\bibliographystyle{siamplain}

\end{document}